# SELF-NORMALIZED PROCESSES: EXPONENTIAL INEQUALITIES, MOMENT BOUNDS AND ITERATED LOGARITHM LAWS

By Victor H. de la Peña,[1] Michael J. Klass[2] and Tze Leung Lai[3]

*Columbia University, University of California at Berkeley and Stanford University*

Self-normalized processes arise naturally in statistical applications. Being unit free, they are not affected by scale changes. Moreover, self-normalization often eliminates or weakens moment assumptions. In this paper we present several exponential and moment inequalities, particularly those related to laws of the iterated logarithm, for self-normalized random variables including martingales. Tail probability bounds are also derived. For random variables $B_t > 0$ and $A_t$, let $Y_t(\lambda) = \exp\{\lambda A_t - \lambda^2 B_t^2/2\}$. We develop inequalities for the moments of $A_t/B_t$ or $\sup_{t \geq 0} A_t/\{B_t(\log \log B_t)^{1/2}\}$ and variants thereof, when $EY_t(\lambda) \leq 1$ or when $Y_t(\lambda)$ is a supermartingale, for all $\lambda$ belonging to some interval. Our results are valid for a wide class of random processes including continuous martingales with $A_t = M_t$ and $B_t = \sqrt{\langle M \rangle_t}$, and sums of conditionally symmetric variables $d_i$ with $A_t = \sum_{i=1}^t d_i$ and $B_t = \sqrt{\sum_{i=1}^t d_i^2}$. A sharp maximal inequality for conditionally symmetric random variables and for continuous local martingales with values in $\mathbf{R}^m$, $m \geq 1$, is also established. Another development in this paper is a bounded law of the iterated logarithm for general adapted sequences that are centered at certain truncated conditional expectations and self-normalized by the square root of the sum of squares. The key ingredient in this development is a new exponential supermartingale involving $\sum_{i=1}^t d_i$ and $\sum_{i=1}^t d_i^2$. A compact law of the iterated logarithm for self-normalized martingales is also derived in this connection.

**1. Introduction.** A prototypical example of self-normalized random variables is Student's $t$-statistic which replaces the population standard devi-

Received November 2002; revised June 2003.
[1]Supported in part by NSF Grants DMS-99-72237 and DMS-02-05791.
[2]Supported in part by NSF Grants DMS-99-72417 and DMS-02-05054.
[3]Supported in part by NSF Grant DMS-00-72523 and NSA Grant MDA 904-00-1-0018.

ation $\sigma$ in the standardized sample mean $\sqrt{n}(\bar{X}_n - \mu)/\sigma$ by the sample standard deviation. More generally, a self-normalized process is of the form $A_t/B_t$, in which $B_t$ is a random variable that estimates some dispersion measure of $A_t$. An important aspect of the theory of self-normalized processes is that we can often dispense with the moment conditions that are needed if $A_t$ is normalized by nonrandom $b_t$ instead, as evidenced by Shao's (1997) large deviation theory for self-normalized sums of i.i.d. random variables without moment conditions. The problem of moment inequalities for self-normalized processes was suggested to the first author in 1990 by J. L. Doob, who pointed out that a key open problem in martingale theory was the development of inequalities for martingales that are analogous to known results in harmonic analysis [see Bañuelos and Moore (1999) for results in this direction].

In recent years, there has been increasing interest in limit theorems and moment bounds for self-normalized sums of i.i.d. zero-mean random variables $X_i$. In particular, Bentkus and Götze (1996) derive a Berry–Esseen bound for Student's $t$-statistic, and Giné, Götze and Mason (1997) prove that the $t$-statistic has a limiting standard normal distribution if and only if $X_1$ is in the domain of attraction of a normal law, by making use of exponential and $L_p$-bounds for the self-normalized sums $U_n = S_n/V_n$, where $S_n = \sum_{i=1}^n X_i$ and $V_n^2 = \sum_{i=1}^n X_i^2$. Egorov (1998) gives exponential inequalities for a centered variant of $U_n$. To see the connection between the $t$-statistic $T_n$ and the self-normalized sum $U_n$, observe that

$$(1.1) \qquad T_n = \frac{S_n/V_n}{\sqrt{\{n - (S_n/V_n)^2\}/(n-1)}}.$$

A recent paper of Caballero, Fernandez and Nualart (1998) contains moment inequalities for a continuous martingale over its quadratic variation and uses these results to show that if $\{M_t, t \geq 0\}$ is a continuous martingale null at zero, then for each $1 \leq p < q$, there exists a universal constant $C = C(p,q)$ such that

$$(1.2) \qquad \left\| \frac{M_t}{\langle M \rangle_t} \right\|_p \leq C \left\| \frac{1}{\langle M \rangle_t^{1/2}} \right\|_q.$$

Related work in Revuz and Yor [(1999), page 168] for continuous local martingales establishes for all $p > q > 0$ the existence of a constant $C_{pq}$ such that

$$(1.3) \qquad E \frac{(\sup_{s<\infty} |M_s|)^p}{\langle M \rangle_\infty^{q/2}} \leq C_{pq} E \left( \sup_{s<\infty} |M_s| \right)^{p-q}.$$

It is important to point out that neither (1.2) nor (1.3) provide bounds for what is arguably the most important case of inequalities of this type, namely $p = q$. Bounds on $E(|M_t|^p/\langle M\rangle_t^{p/2})$ are of particular interest because of their connection with the central limit theorem, as noted earlier in the case of self-normalized sums of i.i.d. random variables. For discrete-time martingales $\{\sum_{i=1}^n d_i, \mathcal{F}_n, n \geq 1\}$, de la Peña (1999) provides exponential bounds for the tail probabilities of $\sum_{i=1}^n d_i/(\alpha + \beta V_n^2)$, where $V_n^2 = \sum_{i=1}^n E(d_i^2|\mathcal{F}_{i-1})$ and $\beta > 0$, $\alpha \geq 0$. In view of the law of the iterated logarithm (LIL), it is of interest to use $V_n$ or $V_n\sqrt{2\log\log V_n}$ (instead of $V_n^2$) to self-normalize $\sum_{i=1}^n d_i$.

Motivated by these developments, we establish in this paper analogous exponential and $L_p$-bounds for a martingale divided by the square root of its quadratic variation or its conditional variance. We start by considering two random variables $A$ and $B$ with $B > 0$ such that

$$(1.4) \qquad E\exp\left\{\lambda A - \frac{\lambda^2}{2}B^2\right\} \leq 1 \qquad \text{for all } \lambda \in \mathbf{R}.$$

Note that if we were allowed to maximize over $\lambda$ inside the expectation, then the maximizing value $\lambda = A/B^2$ would give us $E\exp(A^2/2B^2) \leq 1$, which in turn would imply that $P(A/B \geq x) \leq \exp(-x^2/2)$. Although we cannot interchange the order of $\max_\lambda$ and $E$, we can integrate over $\lambda$ with respect to a probability measure $F$ and interchange the order of integration with respect to $P$ and $F$. This approach is used in Section 2 to derive not only tail probability bounds for $A/B$ but also $L_p$ and exponential bounds for $A/\sqrt{B^2 + (EB)^2}$, and in Section 3 to obtain iterated logarithm bounds for the moments of $A^+/B$. Section 3 further extends the results to the case where (1.4) is replaced by

$$(1.5) \qquad E\exp\{\lambda A - \Phi(\lambda B)\} \leq c \qquad \text{for all } 0 < \lambda < \lambda_0,$$

in which $\Phi$ is assumed to be any nonnegative, strictly convex function on $[0,\infty)$ such that $\Phi(0) = 0$, $\lim_{x\to\infty}\Phi(x) = \infty$ and $\limsup_{x\to\infty}\Phi''(x) < \infty$. Important special cases of such $\Phi$ are $\Phi_r(x) = x^r/r$ with $1 < r \leq 2$.

We next replace the random variables $A$ and $B$ by random processes $A_t$ and $B_t$ and, accordingly, replace (1.5) by

(1.6) $\{\exp(\lambda A_t - \Phi_r(\lambda B_t)), \ t \in T\}$ is a supermartingale for all $0 < \lambda < \lambda_0$,

in which $T$ is either $\{0,1,2,\ldots\}$ or $[0,\infty)$. Section 4 proves an expectation form of the LIL (Theorem 4.1) and develops maximal inequalities under this assumption. Moreover, the case $r = 2$ and $\lambda_0 = \infty$ in (1.6) with "supermartingale" replaced by "martingale" yields a formula for certain boundary crossing probabilities of continuous local martingales taking values in $\mathbf{R}^m$, as shown in Corollary 4.3. Motivated by the LILs for self-normalized sums of certain classes of i.i.d. random variables due to Griffin and Kuelbs (1989, 1991), Shao (1997) and Gine and Mason (1998) and extensions by Jing,

Shao and Wang (2003) to sums of independent zero-mean random variables satisfying a Lindeberg-type condition, we study almost sure LILs for self-normalized (discrete-time) processes in Sections 5 and 6.

When a partial sum of random variables $X_1, X_2, \ldots$ is centered and normalized by a sequence of constants, only under rather special conditions does the usual LIL hold even if the variables are i.i.d. In contrast, we show in Section 5 that there is a universal upper bound of LIL type for the almost sure rate at which such sums can grow after centering by a sum of conditional expectations of suitably truncated variables and normalizing by the square root of the sum of squares of the $X_j$'s. Specifically, let $S_n = X_1 + \cdots + X_n$ and $V_n^2 = X_1^2 + \cdots + X_n^2$, where $\{X_i\}$ is adapted to an increasing sequence $\{\mathcal{F}_i\}$ of $\sigma$-fields. In Section 5 we prove that given any $\lambda > 0$, there exist positive constants $a_\lambda$ and $b_\lambda$ such that $\lim_{\lambda \to 0} b_\lambda = \sqrt{2}$ and

$$(1.7) \quad \limsup_{n \to \infty} \left\{ S_n - \sum_{i=1}^n \mu_i(-\lambda v_n, a_\lambda v_n) \right\} \Big/ \{V_n (\log \log V_n)^{1/2}\} \leq b_\lambda \quad \text{a.s.}$$

on $\{\lim V_n = \infty\}$, where $v_n = V_n (\log \log V_n)^{-1/2}$ and $\mu_i(c, d) = E\{X_i \mathbb{1}(c \leq X_i < d) | \mathcal{F}_{i-1}\}$ for $c < d$. Note that (1.7) is "universal" in the sense that it is applicable to *any* adapted sequence $\{X_i\}$. In particular, suppose $\{S_n, \mathcal{F}_n, n \geq 1\}$ is a supermartingale such that $X_n \geq -m_n$ a.s. for some $\mathcal{F}_{n-1}$-measurable random variable $m_n$ satisfying $P\{0 \leq m_n \leq \lambda v_n \text{ for all large } n\} = 1$. Then (1.7) yields

$$(1.8) \quad \limsup S_n / \{V_n (\log \log V_n)^{1/2}\} \leq b_\lambda \quad \text{a.s. on } \{\lim V_n = \infty\}.$$

We derive in Section 6 the lower half counterpart of (1.8) for the case where $\{S_n, \mathcal{F}_n, n \geq 1\}$ is a martingale such that $|X_n| \leq m_n$ a.s. for some $\mathcal{F}_{n-1}$-measurable $m_n$ with $v_n \to \infty$ and $m_n/v_n \to 0$ a.s. Combining this with (1.8) (with $\lim_{\lambda \to 0} b_\lambda = \sqrt{2}$) then yields

$$(1.9) \quad \limsup S_n / \{V_n (\log \log V_n)^{1/2}\} = \sqrt{2} \quad \text{a.s.}$$

We end this section with various lemmas identifying a large class of random variables satisfying (1.4), (1.5) or (1.6).

LEMMA 1.1. *Let $W_t$ be a standard Brownian motion. Assume that $T$ is a stopping time such that $T < \infty$ a.s. Then $E \exp\{\lambda W_T - \lambda^2 T/2\} \leq 1$ for all $\lambda \in \mathbf{R}$.*

LEMMA 1.2. *Let $M_t$ be a continuous, square-integrable martingale, with $M_0 = 0$. Then*

$$(1.10) \quad \exp\{\lambda M_t - \lambda^2 \langle M \rangle_t / 2\}, \ t \geq 0, \text{ is a supermartingale for all } \lambda \in \mathbf{R}.$$

*If $M_t$ is only assumed to be a continuous local martingale, then (1.10) is also valid (by application of Fatou's lemma).*

LEMMA 1.3. *Let $\{M_t : t \geq 0\}$ be a locally square-integrable martingale, with $M_0 = 0$. Let $\{V_t\}$ be an increasing process, which is adapted, purely discontinuous and locally integrable; let $V^{(p)}$ be its dual predictable projection. Set $X_t = M_t + V_t$, $C_t = \sum_{s \leq t}((\Delta X_s)^+)^2$, $D_t = \{\sum_{s \leq t}((\Delta X_s)^-)^2\}_t^{(p)}$, $H_t = \langle M \rangle_t^c + C_t + D_t$. Then $\exp\{X_t - V_t^{(p)} - \frac{1}{2}H_t\}$ is a supermartingale and, hence,*

(1.11) $\qquad E \exp\{\lambda(X_t - V_t^{(p)}) - \lambda^2 H_t / 2\} \leq 1 \qquad \text{for all } \lambda \in \mathbf{R}.$

Lemma 1.3 is taken from Proposition 4.2.1 of Barlow, Jacka and Yor (1986). A related bound can be found in Lemma 1.5, due to Stout (1973), in which $A_t$ is a discrete-time martingale with bounded increments and $B_t^2$ is a multiple of its conditional variance; see also Kubilius and Mémin (1994). The following lemma holds without any integrability conditions on the variables involved. It is a generalization of the fact that if $X$ is any symmetric random variable, then $A = X$ and $B = X^2$ satisfy condition (1.4). It has a long history, including Wang (1989) and Hitczenko (1990). Hitczenko (1990) proved it for conditionally symmetric martingale difference sequences, and de la Peña (1999) pointed out that the same result still holds without the martingale difference assumption and, hence, without any integrability assumptions.

LEMMA 1.4. *Let $\{d_i\}$ be a sequence of variables adapted to an increasing sequence of $\sigma$-fields $\{\mathcal{F}_i\}$. Assume that the $d_i$'s are conditionally symmetric [i.e., $\mathcal{L}(d_i | \mathcal{F}_{i-1}) = \mathcal{L}(-d_i | \mathcal{F}_{i-1})$]. Then $\exp\{\lambda \sum_{i=1}^n d_i - \lambda^2 \sum_{i=1}^n d_i^2 / 2\}$, $n \geq 1$, is a supermartingale with mean $\leq 1$, for all $\lambda \in \mathbf{R}$.*

Note that any sequence of real-valued random variables $X_i$ can be "symmetrized" to produce an exponential supermartingale satisfying (1.8) by introducing random variables $X_i'$ such that

$$\mathcal{L}(X_n' | X_1, X_1', \ldots, X_{n-1}, X_{n-1}', X_n) = \mathcal{L}(X_n | X_1, \ldots, X_{n-1})$$

and setting $d_n = X_n - X_n'$; see Section 6.1 of de la Peña and Giné (1999). The next two lemmas are related to (1.6).

LEMMA 1.5. *Let $\{d_n\}$ be a sequence of random variables adapted to an increasing sequence of $\sigma$-fields $\{\mathcal{F}_n\}$ such that $E(d_n | \mathcal{F}_{n-1}) \leq 0$ and $d_n \leq M$ a.s. for all $n$ and some nonrandom positive constant $M$. Let $0 < \lambda_0 \leq M^{-1}$, $A_n = \sum_{i=1}^n d_i$, $B_n^2 = (1 + \frac{1}{2}\lambda_0 M) \sum_{i=1}^n E(d_i^2 | \mathcal{F}_{i-1})$, $A_0 = B_0 = 0$. Then $\{\exp(\lambda A_n - \frac{1}{2}\lambda^2 B_n^2), \mathcal{F}_n, n \geq 0\}$ is a supermartingale for every $0 \leq \lambda \leq \lambda_0$.*

LEMMA 1.6. *Let $\{d_n\}$ be a sequence of random variables adapted to an increasing sequence of $\sigma$-fields $\{\mathcal{F}_n\}$ such that $E(d_n|\mathcal{F}_{n-1}) = 0$ and $\sigma_n^2 = E(d_n^2|\mathcal{F}_{n-1}) < \infty$. Assume that there exists a positive constant $M$ such that $E(|d_n|^k|\mathcal{F}_{n-1}) \le (k!/2)\sigma_n^2 M^{k-2}$ a.s. or $P(|d_n| \le M|\mathcal{F}_{n-1}) = 1$ a.s. for all $n \ge 1$, $k > 2$. Let $A_n = \sum_{i=1}^n d_i$, $V_n^2 = \sum_{i=1}^n E(d_i^2|\mathcal{F}_{i-1})$, $A_0 = V_0 = 0$. Then $\{\exp(\lambda A_n - \frac{1}{2(1-M\lambda)}\lambda^2 V_n^2), \mathcal{F}_n, n \ge 0\}$ is a supermartingale for every $0 \le \lambda \le 1/M$.*

Fix any $0 < \rho < 1$. Then Lemma 1.6 implies that (1.4) holds with $A = A_n$ and $B = V_n/\sqrt{\rho}$ for every $0 \le \lambda \le (1-\rho)/M$. The supermartingale in Lemma 1.6 is closely related to martingale extensions of the classical inequalities of Bernstein and Bennett; see Section 8.3 of de la Peña and Giné (1999) for a unified approach to developing such inequalities from corresponding results for sums of independent random variables via decoupling.

**2. Some exponential inequalities.** In this section we present a simple method to derive exponential and $L_p$-bounds for $A/\sqrt{B^2 + (EB)^2}$ under assumption (1.4).

THEOREM 2.1. *Let $B \ge 0$ and $A$ be two random variables satisfying (1.4) for all $\lambda \in \mathbf{R}$. Then for all $y > 0$,*

$$(2.1) \qquad E\frac{y}{\sqrt{B^2+y^2}}\exp\left\{\frac{A^2}{2(B^2+y^2)}\right\} \le 1.$$

*Consequently, if $EB > 0$, then $E\exp(A^2/[4(B^2+(EB)^2)]) \le \sqrt{2}$ and*

$$(2.2) \qquad E\exp(xA/\sqrt{B^2+(EB)^2}) \le \sqrt{2}\exp(x^2) \qquad \text{for all } x > 0.$$

*Moreover, for all $p > 0$,*

$$(2.3) \qquad E(|A|/\sqrt{B^2+(EB)^2})^p \le 2^{p-1/2}p\Gamma(p/2).$$

PROOF. Multiplying both sides of (1.4) by $(2\pi)^{-1/2}y\exp(-\lambda^2 y^2/2)$ (with $y > 0$) and integrating over $\lambda$, we obtain by using Fubini's theorem that

$$1 \ge \int_{-\infty}^{\infty} E\frac{y}{\sqrt{2\pi}}\exp\left(\lambda A - \frac{\lambda^2}{2}B^2\right)\exp\left(-\frac{\lambda^2 y^2}{2}\right)d\lambda$$

$$= E\left[\frac{y}{\sqrt{B^2+y^2}}\exp\left\{\frac{A^2}{2(B^2+y^2)}\right\}\right.$$

$$\left.\times \int_{-\infty}^{\infty} \frac{\sqrt{B^2+y^2}}{\sqrt{2\pi}}\exp\left\{-\frac{B^2+y^2}{2}\left(\lambda^2 - 2\frac{A}{B^2+y^2}\lambda + \frac{A^2}{(B^2+y^2)^2}\right)\right\}d\lambda\right]$$

$$= E\left[\frac{y}{\sqrt{B^2+y^2}}\exp\left\{\frac{A^2}{2(B^2+y^2)}\right\}\right],$$

proving (2.1). By Schwarz's inequality and (2.1),

$$E\exp\left\{\frac{A^2}{4(B^2+y^2)}\right\}$$
$$\leq \left\{\left(E\frac{y\exp\{A^2/(2(B^2+y^2))\}}{\sqrt{B^2+y^2}}\right)\left(E\sqrt{\frac{B^2+y^2}{y^2}}\right)\right\}^{1/2}$$
$$\leq \left(E\sqrt{\frac{B^2}{y^2}+1}\right)^{1/2} \leq \left(E\left(\frac{B}{y}+1\right)\right)^{1/2} \leq \sqrt{2} \quad \text{for } y = EB.$$

To prove (2.2) and (2.3), we assume without loss of generality that $EB < \infty$. Using the inequality $|ab| \leq \frac{a^2+b^2}{2}$, with $a = \sqrt{2c}A/\sqrt{B^2+(EB)^2}$ and $b = x/\sqrt{2c}$, we get $xA/\sqrt{B^2+(EB)^2} \leq \frac{cA^2}{B^2+(EB)^2} + \frac{x^2}{4c}$, which in the case $c = 1/4$ yields

$$E\exp\left\{\frac{xA}{\sqrt{B^2+(EB)^2}}\right\} \leq E\exp\left\{\frac{cA^2}{B^2+(EB)^2} + \frac{x^2}{4c}\right\} \leq \sqrt{2}\exp(x^2),$$

proving (2.2). Moreover, by Markov's inequality, $P(|A|/\sqrt{B^2+(EB)^2} \geq x) \leq \sqrt{2}\exp(-x^2/4)$ for all $x > 0$. Combining this with the formula $EU^p = \int_0^\infty px^{p-1}P(U>x)\,dx$ for any $U \geq 0$, we obtain

$$E(|A|/\sqrt{B^2+(EB)^2})^p \leq \sqrt{2}\int_0^\infty px^{p-1}\exp(-x^2/4)\,dx = 2^{p-1/2}p\Gamma(p/2).$$

$\square$

Another application of the basic inequality (2.1) is the following.

COROLLARY 2.2. *Let $B \geq 0$ and $A$ be two random variables satisfying (1.4) for all $\lambda \in \mathbf{R}$. Then for all $x \geq \sqrt{2}$, $y > 0$ and $p > 0$,*

$$(2.4) \quad P\left(|A|/\sqrt{(B^2+y)\left(1+\frac{1}{2}\log\left(\frac{B^2}{y}+1\right)\right)} \geq x\right) \leq \exp\left(-\frac{x^2}{2}\right),$$

$$(2.5) \quad E\left(|A|/\sqrt{(B^2+y)\left(1+\frac{1}{2}\log\left(\frac{B^2}{y}+1\right)\right)}\right)^p \leq 2^{p/2} + 2^{(p-2)/2}p\Gamma\left(\frac{p}{2}\right).$$

PROOF. Note that for $x \geq \sqrt{2}$ and $y > 0$,

$$P\left\{\frac{A^2}{2(B^2+y)} \geq \frac{x^2}{2}\left(1+\frac{1}{2}\log\left(\frac{B^2}{y}+1\right)\right)\right\}$$

$$\leq P\left\{\frac{A^2}{2(B^2+y)} \geq \frac{x^2}{2} + \frac{1}{2}\log\left(\frac{B^2}{y}+1\right)\right\}$$

$$\leq \exp\left(-\frac{x^2}{2}\right) E \frac{\sqrt{y}\exp\{A^2/(2(B^2+y))\}}{\sqrt{B^2+y}} \leq \exp\left(-\frac{x^2}{2}\right),$$

in which the last inequality follows from (2.1). The proof of (2.5) makes use of (2.4) and is similar to that of (2.3). □

**3. Iterated logarithm bounds for moments of self-normalized variables and their generalizations.** In this section we present bounds for $Eh(A^+/B)$ in terms of $E\{H(B)\}$, where $H$ is a function that depends on $h$. The basic results are Theorems 3.3 and 3.6. Applications of these results are given in Examples 3.4, 3.5 and 3.8, which relate, in particular, the $p$th absolute moment of $A^+/B$ to that of the iterated logarithm $\sqrt{\log\log(B \vee B^{-1} \vee e^2)}$. A variant of Theorem 3.3 has been derived by a different argument in Theorem 1 of de la Peña, Klass and Lai (2000) and Lemmas 3.1 and 3.2 below provide the proofs of Lemmas 2 and 3 of that paper. The main objective of this section is to develop an analogous result that requires (1.4) to hold only for the restricted range $0 < \lambda < \lambda_0$, thereby widely expanding the applicability of our approach. In particular, this extension (given in Theorem 3.6), together with Lemma 1.5, provides moment bounds for a wide class of discrete-time martingales self-normalized by the square root of the conditional variance, thereby connecting our results to LILs. Stout (1973) and Einmahl and Mason (1989) have used this type of self-normalization for LILs of martingales.

Let $L:(0,\infty) \to (0,\infty)$ be a nondecreasing function such that

(3.1) $\qquad L(cy) \leq 3cL(y) \qquad$ for all $c \geq 1$ and $y > 0$,

(3.2) $\qquad L(y^2) \leq 3L(y) \qquad$ for all $y \geq 1$,

(3.3) $\qquad \int_1^\infty \frac{dx}{xL(x)} = \frac{1}{2}.$

An example satisfying (3.1)–(3.3) is the function

(3.4) $\quad L(y) = \beta\{\log(y+\alpha)\}\{\log\log(y+\alpha)\}\{\log\log\log(y+\alpha)\}^{1+\delta},$

where $\delta > 0$, $\alpha$ is chosen sufficiently large to ensure (3.1), (3.2) and $\beta$ is a normalizing constant so that (3.3) holds.

LEMMA 3.1. *Let $\gamma \geq 1$. Then $yL(y/B \vee B/y) \leq 3\gamma\{L(\gamma) \vee L(B \vee B^{-1})\}$ for any $0 < y \leq \gamma$ and $B > 0$. Consequently, for any $A \geq B > 0$ and any $-\frac{A}{B} < x \leq 0$,*

$$(3.5)\ \left(x+\frac{A}{B}\right)L\left(\frac{x+A/B}{B} \vee \frac{B}{x+A/B}\right) \leq 3\frac{A}{B}\left\{L\left(\frac{A}{B}\right) \vee L\left(B \vee \frac{1}{B}\right)\right\}.$$

PROOF. First consider the case $y \leq 1$. From (3.1) and the fact that $L$ is nondecreasing, it follows that

$$yL\left(\frac{y}{B} \vee \frac{B}{y}\right) \leq yL\left(\frac{1}{y}\left(\frac{1}{B} \vee B\right)\right) \leq 3L\left(B \vee \frac{1}{B}\right).$$

For the remaining case $1 < y \leq \gamma$, since $L$ is nondecreasing, we have

$$yL\left(\frac{y}{B} \vee \frac{B}{y}\right) \leq \gamma L\left(\gamma\left(\frac{1}{B} \vee B\right)\right)$$
(3.6)
$$\leq \gamma\left\{L(\gamma^2) \vee L\left(\left(B \vee \frac{1}{B}\right)^2\right)\right\} \leq 3\gamma\left\{L(\gamma) \vee L\left(B \vee \frac{1}{B}\right)\right\},$$

where the last inequality follows from (3.2). □

LEMMA 3.2. *Let $B > 0$ and $A$ be random variables satisfying (1.4) for all $\lambda > 0$. Define*

(3.7) $$g(x) = \frac{\exp\{x^2/2\}}{x}\mathbb{1}(x \geq 1).$$

*Then*

$$E\frac{g(A/B)}{L(A/B) \vee L(B \vee 1/B)} \leq \frac{3}{\int_0^1 \exp(-x^2/2)\,dx}.$$

PROOF. By a change of variables, $\int_0^1 (\lambda L(1/\lambda))^{-1}\,d\lambda = \int_1^\infty (\lambda L(\lambda))^{-1}\,d\lambda = \frac{1}{2}$. Let

(3.8) $$f(\lambda) = \frac{1}{\lambda L(\max\{\lambda, 1/\lambda\})}, \qquad \lambda > 0.$$

Then $\int_0^\infty f(\lambda)\,d\lambda = \int_0^1 f(\lambda)\,d\lambda + \int_1^\infty f(\lambda)\,d\lambda = 1$, so $f$ is a density function on $(0, \infty)$. Therefore, integrating (1.4) with respect to this probability measure yields

$$1 \geq E\int_0^\infty \frac{\exp\{Ax - (B^2x^2/2)\}}{xL(x \vee 1/x)}\,dx$$

$$= E\int_0^\infty \frac{\exp\{Ay/B - (y^2/2)\}}{yL(y/B \vee B/y)}\,dy \qquad \text{(letting } y = Bx\text{)}$$

$$\geq E\left\{\exp\left(\frac{A^2}{2B^2}\right)\right\}$$

$$\times \int_{-A/B}^\infty \frac{\exp\{-(x^2/2)\}}{(x + A/B)L(\{(x + A/B)/B\} \vee \{B(x + A/B)\})}\mathbb{1}\left(\frac{A}{B} \geq 1\right)dx$$

$$\left(\text{letting } x = y - \frac{A}{B}\right)$$

$$\geq E\left\{\exp\left(\frac{A^2}{2B^2}\right)\right\} \int_{-1}^{0} \frac{\exp\{-(x^2/2)\}\,dx}{3(A/B)(L(A/B) \vee L(B \vee 1/B))} \mathbb{1}\left(\frac{A}{B} \geq 1\right)$$

[by (3.5)]

$$= \left\{\frac{1}{3}\int_0^1 \exp\left(-\frac{x^2}{2}\right) dx\right\} E \frac{g(A/B)}{L(A/B) \vee L(B \vee 1/B)}. \qquad \square$$

We next derive a bound on $Eh(A^+/B)$ by making use of Lemma 3.2 for nondecreasing functions $h$ that do not grow faster than $g/L$.

THEOREM 3.3. *Let $L:(0,\infty) \to (0,\infty)$ be a nondecreasing function satisfying (3.1)–(3.3). Define $g$ by (3.7). Let $h$ be a nondecreasing function on $[0,\infty)$ such that for some $x_0 \geq 1$ and $c > 0$,*

(3.9) $\qquad\qquad 0 < h(x) \leq cg(x)/L(x) \qquad$ *for all $x \geq x_0$.*

*Let $q$ be a strictly increasing, continuous function on $[0,\infty)$ such that for some $\bar{c} \geq c$,*

(3.10) $\qquad\qquad L(x) \leq q(x) \leq \dfrac{\bar{c}g(x)}{h(x)} \qquad$ *for all $x \geq x_0$.*

*Let $B > 0$ and $A$ be random variables satisfying (1.4) for all $\lambda > 0$. Then*

(3.11) $\qquad Eh(A^+/B) \leq 4\bar{c} + h(x_0) + Eh(q^{-1}(L(B \vee B^{-1}))).$

*Consequently, $Eh(A^+/B) < \infty$ if $Eh(q^{-1}(L(B \vee B^{-1}))) < \infty$.*

PROOF. By Lemma 3.2,

$$E\frac{g(A^+/B)}{L(A/B) \vee L(B \vee 1/B)} \leq 4.$$

Let $Q = \{L(B \vee \frac{1}{B}) \leq q(\frac{A}{B})\}$. Then, $Eh(A^+/B)$ is majorized by

$$h(x_0) + E\frac{h(A^+/B)\mathbb{1}(Q)\mathbb{1}(A/B \geq x_0)}{g(A/B)/(L(A/B) \vee L(B \vee 1/B))}$$

$$\times \left(\frac{g(A/B)}{L(A/B) \vee L(B \vee 1/B)}\right) + Eh\left(\frac{A^+}{B}\right)\mathbb{1}(Q^c)\mathbb{1}\left(\frac{A}{B} \geq x_0\right)$$

$$\leq h(x_0) + \sup_{y \geq x_0} \frac{h(y)(L(y) \vee q(y))}{g(y)}$$

$$\times E\left(\frac{g(A/B)}{L(A/B) \vee L(B \vee 1/B)}\right) + Eh\left(q^{-1}\left(L\left(B \vee \frac{1}{B}\right)\right)\right)$$

$$\leq h(x_0) + 4\sup_{y \geq x_0} \frac{h(y)q(y)}{g(y)} + Eh\left(q^{-1}\left(L\left(B \vee \frac{1}{B}\right)\right)\right).$$

□

To apply Theorem 3.3, one can take $L$ as given by (3.4) and choose $q^{-1}$ that grows as slowly as possible (or equivalently, $q$ that grows as rapidly as possible) subject to the constraint (3.10).

EXAMPLE 3.4. Define $L$ by (3.4) and let $h(x) = x^p$ for $x \geq 0$, with $p > 0$. Then (3.9) clearly holds with $c = 1$ and $x_0$ sufficiently large, for which (3.10) also holds with $q(x) = g(x)/h(x) = \exp(x^2/2)/x^{p+1}$. In this case,

$$q^{-1}(y) = \{2\log y + (p+1+o(1))\log\log y\}^{1/2} \quad \text{as } y \to \infty.$$

Since $L(x) \sim \beta(\log x)(\log\log x)(\log\log\log x)^{1+\delta}$ as $x \to \infty$, Theorem 3.3 yields

(3.12) $\quad E(A^+/B)^p < \infty \quad \text{if } E\{\log(|\log(B \vee B^{-1})| \vee e)\}^{p/2} < \infty,$

for random variables $B > 0$ and $A$ satisfying (1.4) for all $\lambda > 0$.

EXAMPLE 3.5. Let $0 < \theta < 1$ and $h(x) = \exp(\theta x^2/2)$ for $x \geq 0$. Define $L$ by (3.4). Then (3.9) holds with $c = 1$ and $x_0$ sufficiently large, for which (3.10) also holds with $q(x) = g(x)/h(x) = x^{-1}\exp\{(1-\theta)x^2/2\}$. In this case, $h(q^{-1}(y)) = O(\{y(\log y)^{1/2}\}^{\theta/(1-\theta)})$. Therefore, if $B > 0$ and (1.4) holds for all $\lambda > 0$, then by Theorem 3.3,

(3.13)
$$E\exp\left(\frac{\theta}{2}\left(\frac{A^+}{B}\right)^2\right) < \infty$$
$$\text{if } E\{(\log \tilde{B})(\log\log \tilde{B})^{3/2}(\log\log\log \tilde{B})^{1+\delta}\}^{\theta/(1-\theta)} < \infty$$

for some $\delta > 0$, where $\tilde{B} = B \vee B^{-1} \vee e^3$.

The following theorem modifies and broadly extends Theorem 3.3 by requiring (1.4) to hold only for the restricted range $0 \leq \lambda \leq \lambda_0$. An example where this appears naturally can be found in Lemma 1.5, where $A$ is a martingale and $B^2$ is a multiple of its conditional variance. Theorem 3.6 also generalizes (1.4) by replacing the quadratic function $\lambda^2 B^2/2$ and the upper bound 1 in (1.4) by a convex function $\Phi(\lambda B)$ and a finite positive constant $c$. Unlike Theorem 3.3 that involves a single function $q$ to give the bound (3.11), Theorem 3.6 uses a family of functions $q_b$. The wider range of applications that will be explored in Section 4 justifies the additional technical work required for the theorem. The proof employs different analyses of $A/B$ for small and large $B$, incorporating a Taylor expansion of $\Phi$ for large $B$. In addition, as before, Fubini's theorem allows us to treat the random variables involved as though they were constants.

THEOREM 3.6. *Suppose that $\Phi(\cdot)$ is a continuous function with $\Phi'(x)$ strictly increasing, continuous and positive for $x > 0$, with $\lim_{x \to \infty} \Phi(x) = \infty$ and $\sup_{x>0} \Phi''(x) < \infty$. Suppose $B > 0$ and $A$ are random variables such that there exists $c > 0$ for which*

$$(3.14) \qquad E\exp\{\lambda A - \Phi(\lambda B)\} \leq c \qquad \text{for all } 0 < \lambda < \lambda_0.$$

*For $w > \Phi'(1)$, define $y_w$ by the equation $\Phi'(y_w) = w$, and let*

$$(3.15) \qquad g_\Phi(w) = y_w^{-1} \exp\{wy_w - \Phi(y_w)\}.$$

*Let $\eta > \tilde\eta > 0$. Let $h: [0, \infty) \to (0, \infty)$ be a nondecreasing function. For $b \geq \eta$, let $q_b$ be a strictly increasing, continuous function on $(0, \infty)$ such that for some $\tilde c > 0$ and $w_0 > \Phi'(2)$,*

$$(3.16) \qquad q_b(w) \leq \tilde c\{g_\Phi(w)\mathbb{1}(y_w \leq \lambda_0 b) + e^{\lambda_0 \tilde\eta w}\mathbb{1}(y_w > \lambda_0 b)\}/h(w)$$

*for all $w \geq w_0$.*

*Let $L: (0, \infty) \to (0, \infty)$ be a nondecreasing function satisfying (3.1)–(3.3). Then there exists a constant $C$ depending only on $\lambda_0, \eta, \tilde\eta, c, \tilde c$ and $\Phi$ such that*

$$(3.17) \qquad Eh(A^+/(B \vee \eta)) \leq C + h(w_0) + Eh(q_{B \vee \eta}^{-1}(L(B \vee \eta))).$$

PROOF. Note that Lemma 3.2 transforms the inequality constraints (1.4) for all $\lambda > 0$ into a single expectation inequality primarily involving a rather rapidly growing function of $A/B$ and a slowly varying function $L$ of $B \vee \frac{1}{B}$. This result is employed in the proof of Theorem 3.3 to bound a quantity of the form $Eh(A^+/B)$ by a constant plus $Eh(q^{-1}(L(B \vee \frac{1}{B})))$. To duplicate this approach when (1.4) holds only for $0 < \lambda < \lambda_0$, we first derive an analog of Lemma 3.2 by splitting $A/B \geq w_0$ into two cases: $y_{A/B} > \lambda_0 B$ and $y_{A/B} \leq \lambda_0 B$. Moreover, we need to replace $B$ by $B \vee \eta$. Since $\Phi(x)$ is increasing in $x > 0$, (3.14) also holds with $B$ replaced by $B \vee \eta$ and, therefore, we shall assume without loss of generality that $B \geq \eta$. Integrating (3.14) with respect to the probability measure defined by the density function (3.8) yields

$$(3.18) \quad c \geq E \int_0^{\lambda_0} \frac{\exp\{\lambda A - \Phi(\lambda B)\}}{\lambda L(\lambda \vee \lambda^{-1})} d\lambda = E \int_0^{\lambda_0 B} \frac{\exp\{xA/B - \Phi(x)\}}{xL(x/B \vee B/x)} dx.$$

Our first variant of Lemma 3.2, given in (3.19), provides an exponential bound for $A/B$ when $\lambda_0 B < y_{A/B}$. Observe that using the definition of $y_w$, we have that $x\frac{A}{B} - \Phi(x)$ increases in $x$ for $x \leq y_{A/B}$, and decreases in $x$ for $x \geq y_{A/B}$. Take any $0 < \tilde{\eta} < \eta$, and let $\lambda_1 = \lambda_0 \vee \lambda_0^{-1} \vee \tilde{\eta}$. Since $B \geq \eta > \tilde{\eta}$, it follows from (3.18) and (3.1) that

$$c \geq E \int_{\lambda_0 \tilde{\eta}}^{\lambda_0 \eta} \frac{\exp\{xA/B - \Phi(x)\}}{xL(x/B \vee B/x)} dx \mathbb{1}\left(\frac{A}{B} \geq w_0\right) \mathbb{1}(y_{A/B} > \lambda_0 B)$$

$$(3.19) \quad \geq E \int_{\lambda_0 \tilde{\eta}}^{\lambda_0 \eta} \frac{\exp\{\lambda_0 \tilde{\eta} A/B - \Phi(\lambda_0 \eta)\}}{L(\lambda_0 \vee B/(\lambda_0 \tilde{\eta}))} \frac{dx}{x} \mathbb{1}\left(\frac{A}{B} \geq w_0\right) \mathbb{1}(y_{A/B} > \lambda_0 B)$$

$$\geq \frac{e^{-\Phi(\lambda_0 \tilde{\eta})}}{3\lambda_1/\tilde{\eta}} \log\left(\frac{\eta}{\tilde{\eta}}\right) E \frac{e^{\lambda_0 \tilde{\eta} A/B}}{L(B)} \mathbb{1}\left(\frac{A}{B} \geq w_0\right) \mathbb{1}(y_{A/B} > \lambda_0 B).$$

Our second variant of Lemma 3.2, given in (3.21), bounds $A/B$ when $\lambda_0 B \geq y_{A/B}$. Since $w_0 > \Phi'(2), y_{w_0} > 2$. Define

$$(3.20) \qquad a_* = \sup\{a \leq 1 : a^2 \Phi''(x) \leq 1 \text{ for all } x > y_{w_0} - a\}.$$

Note that $a_* > 0$ and $y_{w_0} - a_* > 1$. Since $\Phi'(y_w) - w = 0$, a two-term Taylor expansion for $w \geq w_0$ and $x \in (y_w - a_*, y_w)$ yields

$$wx - \Phi(x) = wy_w - \Phi(y_w) - \frac{(x - y_w)^2}{2} \Phi''(\xi^*)$$

$$\geq wy_w - \Phi(y_w) - \frac{(x - y_w)^2}{2a_*^2},$$

in which $\xi^*$ lies between $x$ and $y_w$. The last inequality follows from (3.16) and (3.20), noting that $\xi^* > x > y_w - a_* \geq y_{w_0} - a_*$. It then follows from (3.18) that

$$c \geq E\bigg[\mathbb{1}\bigg(y_{A/B} \leq \lambda_0 B, \frac{A}{B} \geq w_0\bigg)$$

$$\times \int_{y_{A/B} - a_*}^{y_{A/B}} \frac{\exp\{(A/B)y_{A/B} - \Phi(y_{A/B}) - (x - y_{A/B})^2/(2a_*^2)\}}{xL(x/B \vee B/x)} dx\bigg]$$

$$\geq E\bigg[\mathbb{1}\bigg(y_{A/B} \leq \lambda_0 B, \frac{A}{B} \geq w_0\bigg)$$

$$\times \frac{\exp\{(A/B)y_{A/B} - \Phi(y_{A/B})\}}{y_{A/B}\{L(\lambda_0 \vee B)\}} \int_{y_{A/B} - a_*}^{y_{A/B}} \exp\bigg\{-\frac{(x - y_{A/B})^2}{2a_*^2}\bigg\} dx\bigg],$$

using $x > y_{A/B} - a_* \geq y_{w_0} - a_* > 1$ so that $\frac{B}{x} < B$. From Lemma 3.1 and the fact that $B \geq \eta$, we have $L(\lambda_0 \vee B) \leq 3(1 \vee \frac{\lambda_0}{\eta}) L(B)$. Hence,

$$c \geq E\bigg[\mathbb{1}\bigg(y_{A/B} \leq \lambda_0 B, \frac{A}{B} \geq w_0\bigg)$$

$$(3.21) \quad \times \frac{\exp\{(A/B)y_{A/B} - \Phi(y_{A/B})\}}{3y_{A/B}(1 \vee (\lambda_0/\eta))L(B)} a_* \int_0^1 \exp\left(-\frac{z^2}{2}\right) dz \Bigg]$$

$$\geq \frac{a_*}{4(1 \vee (\lambda_0/\eta))} E \frac{g_\Phi(y_{A/B}) \mathbb{1}(y_{A/B} \leq \lambda_0 B, A/B \geq w_0)}{L(B)}.$$

Let $Q = \{L(B) \leq q_B(A/B)\}$. Then rewriting (3.16) as an upper bound for $h$ and using the definition of $Q$, we can majorize $Eh(A^+/B)$ by

$$h(w_0) + \tilde{c}E\left[\mathbb{1}(Q)\left\{\frac{g_\Phi(A/B)}{L(B)}\mathbb{1}\left(\frac{A}{B} \geq w_0, y_{A/B} \leq \lambda_0 B\right)\right.\right.$$

$$\left.\left. + \frac{e^{\lambda_0 \tilde{\eta} A/B}}{L(B)}\mathbb{1}\left(\frac{A}{B} \geq w_0, y_{A/B} > \lambda_0 B\right)\right\}\right]$$

$$+ Eh\left(\frac{A}{B}\right)\mathbb{1}\left(Q^c \cap \left\{\frac{A}{B} \geq w_0\right\}\right)$$

$$\leq h(w_0) + C + Eh(q_B^{-1}(L(B))),$$

in which the inequality follows from (3.16), (3.21) and (3.19). □

REMARK 3.7. In the case $\lambda_0 = \infty$ [as in Theorem 3.3 for which $\Phi(x) = x^2/2$], the bounds (3.18) and (3.19) are not needed and the result for general $\Phi$ is similar to (3.11) in Theorem 3.3. The main difference between (3.11) and (3.17) lies in $q^{-1}$ in (3.11) versus the more elaborate $q_{B \vee \eta}^{-1}$ in (3.17) to incorporate both (3.19) and (3.21).

The next example is designed to exploit the form of $q_b(w)$ of Theorem 3.6 [see (3.16)].

EXAMPLE 3.8. Lemmas 1.5 and 1.6 give examples of $(A, B)$ satisfying (1.4) only for $0 \leq \lambda \leq \lambda_0$. Thus, (3.14) holds with $\Phi(x) = x^2/2$ and $g_\Phi$ reduces to the function $g$ defined by (3.7) in this case, noting that $y_w = w$. Define $L$ by (3.4). First let $h(x) = x^p$ for $x \geq 0$, with $p > 0$. For $b \geq \eta > \tilde{\eta} > 0$, let $q_b$ be a strictly increasing function on $(0, \infty)$ such that for all large $b$,

$$(3.22) \quad \begin{aligned} q_b(w) &= e^{w^2/2}/w^{p+1} && \text{if } w \leq \lambda_0(\tilde{\eta}b)^{1/2}, \\ &\leq e^{w^2/2}/w^{p+1} && \text{if } \lambda_0(\tilde{\eta}b)^{1/2} < w \leq \lambda_0 b, \\ &= e^{\lambda_0 \tilde{\eta} w}/w^p && \text{if } w > \lambda_0 b. \end{aligned}$$

Then (3.16) holds with $\tilde{c} = 1$. From (3.4) and (3.22), it follows that $q_b^{-1}(L(b)) \sim (2 \log \log b)^{1/2}$ as $b \to \infty$. Therefore, (3.12) still holds with $B$ replaced by $B \vee \eta$ even though (1.4) holds only for $0 \leq \lambda \leq \lambda_0$. Similarly, letting $h(x) = e^{\zeta x}$ with $0 < \zeta < \lambda_0 \tilde{\eta}$, it follows from Theorem 3.6 that

$$(3.23) \quad \begin{aligned} &E \exp(\zeta A^+/(B \vee \eta)) < \infty \\ &\quad \text{if } E \exp\{\zeta[2(\log \log \tilde{B})(\log \log \log \tilde{B})^{1+\delta}]^{1/2}\} \end{aligned}$$

for some $\delta > 0$, where $\tilde{B} = B \vee e^3$. One such choice of $q_b$ that satisfies (3.22) for sufficiently large $b$ is to let $q_b(w) = w^{-p} \exp(f^2(w))$ for $\lambda_0(\tilde{\eta}b)^{1/2} < w \leq \lambda_0 b$, where $f$ is linear on $[\lambda_0(\tilde{\eta}b)^{1/2}, \lambda_0 b]$ and is uniquely determined by requiring $q_b$ to be continuous. In this case, it can be shown that $f^2(w) \leq w^2/2 - \log w$ for $\lambda_0(\tilde{\eta}b)^{1/2} \leq w \leq \lambda_0 b$ if $b$ is sufficiently large, noting that the slope of $f$ is $\{1 + o(1) - 1/\sqrt{2}\}\sqrt{\tilde{\eta}/b}$ and, therefore, $\frac{1}{2}w^2 - \log w - f^2(w)$ is an increasing function of $w \in [\lambda_0(\tilde{\eta}b)^{1/2}, \lambda_0 b]$ for all large $b$.

Another application of Theorem 3.6 involves the more general case of $\Phi(x) = x^r/r$ $(1 < r \leq 2)$, for which

$$(3.24) \quad y_w = w^{1/(r-1)}, \qquad g_\Phi(w) = w^{-1/(r-1)} \exp\{(1 - r^{-1})w^{r/(r-1)}\}.$$

In view of (3.24), it follows from Theorem 3.6, by arguments similar to Example 3.8, that under (3.14) with $\Phi(x) = x^r/r$, we have for any $p > 0$,

$$(3.25) \quad E(A^+/(B \vee \eta))^p < \infty \qquad \text{if } E\{\log^+(\log(B \vee \eta))\}^{p(r-1)/r} < \infty.$$

Moreover, (3.23) still holds if we replace 2 and $1/2$ there by $r/(r-1)$ and its reciprocal, respectively. The following lemma, which provides an analogue of Lemma 1.5 for more general $1 < r \leq 2$ and which self-normalizes $A_n$ by the square root of the square function $\sum_{i=1}^n d_i^2$, gives an exponential supermartingale when the summands $d_i$ of $A_n$ are bounded from below rather than from above.

LEMMA 3.9. *Let $0 < \gamma < 1 < r \leq 2$. Define $c_{\gamma,r} = \max\{c_r, c_r^{(\gamma)}\}$, where*

$$c_r = \inf\{c > 0 : \exp(x - cx^r) \leq 1 + x \text{ for all } x \geq 0\},$$
$$c_r^{(\gamma)} = \inf\{c > 0 : \exp(x - c|x|^r) \leq 1 + x \text{ for all } -\gamma \leq x \leq 0\}.$$

(i) *For all $x \geq -\gamma$, $\exp\{x - c_{\gamma,r}|x|^r\} \leq 1 + x$. Moreover, $c_r \leq (r-1)^{r-1}(2-r)^{2-r}/r$ and*

$$c_r^{(\gamma)} = -\{\gamma + \log(1-\gamma)\}/\gamma^r = \sum_{j=2}^\infty \gamma^{j-r}\Big/j.$$

(ii) *Let $\{d_n\}$ be a sequence of random variables adapted to an increasing sequence of $\sigma$-fields $\{\mathcal{F}_n\}$ such that $E(d_n|\mathcal{F}_{n-1}) \leq 0$ and $d_n \geq -M$ a.s. for all $n$ and some nonrandom positive constant $M$. Let $A_n = \sum_{i=1}^n d_i$, $B_n^r = rc_{\gamma,r} \sum_{i=1}^n |d_i|^r$, $A_0 = B_0 = 0$. Then $\{\exp(\lambda A_n - (\lambda B_n)^r/r), \mathcal{F}_n, n \geq 0\}$ is a supermartingale for every $0 \leq \lambda \leq \gamma M^{-1}$.*

PROOF. The first assertion of (i) follows from the definition of $c_{\gamma,r}$. For $c > 0$, define $g_c(x) = \log(1 + x) - x + c|x|^r$ for $x > -1$. Then $g_c'(x) =$ $|x|^{r-1}\{|x|^{2-r}(1-|x|)^{-1} - cr\}$ for $-1 < x < 0$. Since $|x|^{2-r}/(1-|x|)$ is decreasing in $-1 < x < 0$, $g'_c$ has at most one zero belonging to $(-1, 0)$. Let $c^* = -\{\gamma + \log(1-\gamma)\}/\gamma^r$. Then $g_{c^*}(-\gamma) = 0 = g_{c^*}(0)$. It then follows that $g_{c^*}(x) > 0$ for all $-\gamma < x < 0$ and, therefore, $c^* \geq c_r^{(\gamma)}$. If $c^* > c_r^{(\gamma)}$, then $g_{c_r^{(\gamma)}}(-\gamma) < g_{c^*}(-\gamma) = 0$, contradicting the definition of $c_r^{(\gamma)}$. Hence, $c_r^{(\gamma)} = c^*$. Take any $c \geq (r-1)^{r-1}(2-r)^{2-r}/r$. Then for all $x > 0$,

$$g'_c(x) = \frac{1}{1+x} - 1 + crx^{r-1} \geq \frac{x}{1+x}\left\{-1 + cr\inf_{y>0}(y^{r-2} + y^{r-1})\right\}$$

$$= \frac{x}{1+x}\left\{-1 + \frac{cr}{(r-1)^{r-1}(2-r)^{2-r}}\right\} \geq 0.$$

Since $g_c(0) = 0$, it then follows that $g_c(x) \geq 0$ for all $x \geq 0$. Hence, $c_r \leq (r-1)^{r-1}(2-r)^{2-r}/r$.

To prove (ii), note that since $\lambda d_n \geq -\lambda M \geq -\gamma$ a.s. for $0 \leq \lambda \leq \gamma M^{-1}$, (i) yields

$$E[\exp\{\lambda d_n - c_{\gamma,r}|\lambda d_n|^r\}|\mathcal{F}_{n-1}] \leq E[1 + \lambda d_n|\mathcal{F}_{n-1}] \leq 1 \qquad \text{a.s.} \qquad \square$$

**4. An expectation version of the LIL and maximal inequalities for self-normalized martingales.** In this section we first prove a theorem that provides an expectation form of the upper LIL under the assumption

(4.1)
$$\{\exp(\lambda A_t - \Phi_r(\lambda B_t)), t \in T\}$$

is a supermartingale with mean $\leq 1$ for $0 < \lambda < \lambda_0$,

where $T$ is either $\{0, 1, 2, \dots\}$ (discrete-time case) or $[0, \infty)$ (continuous-time case) and $\Phi_r(x) = x^r/r$ for $1 < r \leq 2$. Applications of the theorem will be given in (4.9)–(4.12). Important special cases of (4.1) have been given in Lemmas 1.5, 1.6 and 3.9. We then develop maximal inequalities for self-normalized processes under (4.1), yielding an almost sure upper LIL in Corollary 4.2 that generalizes a corresponding result of Giné and Mason (1998) for i.i.d. symmetric random variables.

THEOREM 4.1. *Let $T = \{0, 1, 2, \dots\}$ or $T = [0, \infty)$, $1 < r \leq 2$, and $\Phi_r(x) = x^r/r$ for $x > 0$. Let $A_t, B_t$ be stochastic processes (on the same probability space) satisfying (4.1) and such that $B_t$ is positive and nondecreasing in $t > 0$, with $A_0 = 0$. In the case $T = [0, \infty)$, assume furthermore that $A_t$ and $B_t$ are right-continuous. Let $L : [1, \infty) \to (0, \infty)$ be a nondecreasing function satisfying (3.1)–(3.3). Let $\eta > 0$, $\lambda_0\eta > \varepsilon > 0$, and $h : [0, \infty) \to [0, \infty)$ be a nondecreasing function such that $h(x) \leq e^{\varepsilon x}$ for all large $x$. Then there exists a constant $C$ depending only on $\lambda_0, \eta, r, \varepsilon, h$ and $L$ such that*

(4.2) $$Eh\left(\sup_{t \geq 0}\{A_t(B_t \vee \eta)^{-1}[1 \vee \log^+ L(B_t \vee \eta)]^{-(r-1)/r}\}\right) \leq C.$$

PROOF. It suffices to prove (4.2) with $\sup_{t\geq 0}$ replaced by $\sup_{s\geq t\geq 0}$ for every $s > 0$. Given any $s > 0$, there exists a sequence of nonnegative random times $\tau_n \leq s$ (in general, not stopping times) such that

$$\text{(4.3)} \quad \lim_{n\to\infty} \frac{A^+_{\tau_n}}{(B_{\tau_n} \vee \eta)\{1 \vee \log^+ L(B_{\tau_n} \vee \eta)\}^{(r-1)/r}} = \sup_{0\leq t\leq s} \frac{A_t}{(B_t \vee \eta)\{1 \vee \log^+ L(B_t \vee \eta)\}^{(r-1)/r}},$$

since $A_0 = 0$. As in the proof of Theorem 3.6, we shall assume without loss of generality that $B_t \geq \eta$. Take any $q < 1$ such that $q\lambda_0 \eta > \varepsilon$.

It follows from Lemma 1 of Shao (2000) and Fatou's lemma that for any nonnegative supermartingale $\{Y_t, t \in T\}$ (with right-continuous $Y_t$ in the case $T = [0,\infty)$), $E(\sup_{t\in T} Y_t)^q \leq (1-q)^{-1}(EY_0)^q$. Applying this result to (4.1) and noting that $A_0 = 0$, we obtain that for $0 \leq \lambda \leq \lambda_0$,

$$\text{(4.4)} \quad \begin{aligned} (1-q)^{-1} &\geq E\left(\sup_{t\in T}\exp\{\lambda A_t - \Phi_r(\lambda B_t)\}\right)^q \\ &\geq E\exp\{q[\lambda A_{\tau_n} - \Phi_r(\lambda B_{\tau_n})]\} \\ &= E\exp\{q\lambda A_{\tau_n} - \Psi_r(q\lambda B_{\tau_n})\}, \end{aligned}$$

where $\Psi_r(x) = q^{1-r} x^r / r$.

Let $f_r(w) = \exp\{q(1-r^{-1})w^{r/(r-1)}\}$ for $w > 0$. Note that in the notation of Theorem 3.6, $g_{\Psi_r}(w) = y_w^{-1} f_r(w)$ with $y_w = qw^{1/(r-1)}$. Letting $A = A_{\tau_n}$ and $B = B_{\tau_n}$, it follows from (4.4) and (3.3) that

$$(1-q)^{-1} \geq \int_0^{\lambda_0} E\exp\{q\lambda A - \Psi_r(q\lambda B)\}\frac{d\lambda}{\lambda L(\lambda \vee \lambda^{-1})},$$

which in turn yields the following analogues of (3.19) and (3.21), with $\eta > \tilde{\eta}$:

$$\text{(4.5)} \quad (1-q)^{-1} \geq \frac{e^{-\Psi_r(q\lambda_0 \eta)}}{3\lambda_1(\tilde{\eta}^{-1} \vee \tilde{\eta})} \log\left(\frac{\eta}{\tilde{\eta}}\right) E\frac{e^{q\lambda_0 \tilde{\eta} A/B}}{L(B \vee 1)},$$

$$\text{(4.6)} \quad \begin{aligned} (1-q)^{-1} &\geq E\int_0^{q\lambda_0 B} \frac{\exp\{xA/B - \Psi_r(x)\}}{xL(x/(qB) \vee (qB)/x)}\,dx \\ &\geq c(\lambda_0, q, \eta, r) E\frac{f_r(A/B)}{(A/B)^{1/(r-1)}L(B)}\mathbb{1}\left(\left(\frac{A}{B}\right)^{1/(r-1)} \leq \lambda_0 B\right), \end{aligned}$$

with $q\lambda_0 \tilde{\eta} > \varepsilon$, $\lambda_1 = \lambda_0 \vee \lambda_0^{-1}$ and the constant $c(\lambda_0, q, \eta, r)$ depending only on $\lambda_0, q, \eta$ and $r$. For (4.6), recall that $y_w = qw^{1/(r-1)}$ and $g_{\Psi_r}(w) = y_w^{-1} f_r(w)$.

Take any $\delta < 1$ such that $r(1-\delta)/(r-1) > 1$. Since $q\lambda_0 \tilde{\eta} > \varepsilon$, there exists $x_0 > \lambda_0^{r-1} \vee 1$ such that

$$\text{(4.7)} \quad h(x) \leq e^{q\lambda_0 \tilde{\eta} x}/L(x) < f_r^{1-\delta}(x)/x^{1/(r-1)} \quad \text{for all } x \geq x_0,$$

noting that $L(x) \leq 3xL(1)$ by (3.1). Let
$$F = \{f_r^\delta(A/B) \leq L(B) \vee e\}$$
$$= \{A^+/B \leq [(1 \vee \log L(B))/(\delta q(1-r^{-1}))]^{(r-1)/r}\}.$$
Let $k$ be the smallest integer such that $2^k(r-1) \geq 1$. On $\{A/B \geq x_0 \vee (\lambda_0 B)^{r-1}\}$,
$$L(A/B) \geq L(x_0 \vee (\lambda_0 B)^{r-1}) \geq \tfrac{1}{3}(\lambda_0 \wedge 1)^{r-1} L(1 \vee B^{r-1})$$
$$\geq 3^{-(k+1)}(\lambda_0 \wedge 1)^{r-1} L(1 \vee B),$$
where the last two inequalities follow from (3.1) and (3.2), respectively. From (4.7), it then follows that

$$Eh\left(\frac{A^+}{B\{1 \vee \log^+ L(B)\}^{(r-1)/r}}\right)$$
$$\leq h(x_0) + h(1/[\delta q(1-r^{-1})]^{(r-1)/r}) P(F)$$
$$+ E\mathbb{1}\left(F^c \cap \left\{\frac{A}{B} \geq x_0\right\}\right)$$
(4.8)
$$\times \left\{\frac{e^{q\lambda_0 \tilde{\eta} A/B}}{3^{-(k+1)}(\lambda_0 \wedge 1)^{r-1} L(B \vee 1)} \mathbb{1}\left(\left(\frac{A}{B}\right)^{1/(r-1)} \geq \lambda_0 B\right)\right.$$
$$+ \left.\left(\sup_{x \geq x_0} \frac{h(x) x^{1/(r-1)}}{f_r^{1-\delta}(x)}\right) \frac{f_r(A/B)}{(A/B)^{1/(r-1)} L(B)}\right.$$
$$\left.\times \mathbb{1}\left(\left(\frac{A}{B}\right)^{1/(r-1)} < \lambda_0 B\right)\right\},$$

noting that $f_r^\delta(A/B) > L(B)$ and, therefore, $f_r^{1-\delta}(A/B) < f_r(A/B)/L(B)$ on $F^c$. The desired conclusion follows from (4.5), (4.6) and (4.8). □

Consider the case of continuous local martingales $A_t$. We can apply Theorem 4.1 with $r = 2$ and $B_t = \sqrt{\langle A \rangle_t}$, in view of Lemma 1.2. Putting $h(x) = x^p$ in (4.2), with $L(x)$ given by (3.4) in this case, yields the following extension of (1.3) to the case $q = p$: There exists for every $p > 0$ an absolute constant $C_p$ such that

(4.9) $$E\left(\sup_{t \geq 0} \frac{A_t^+}{\{\langle A \rangle_t \log\log(\langle A \rangle_t \vee e^2)\}^{1/2}}\right)^p \leq C_p.$$

Since (4.1) holds for all $\lambda_0 > 0$ by Lemma 1.2, we can, in fact, set $\lambda_0 = \infty$ in (4.8) with $r = 2$ to replace it by

$$Eh\left(\frac{A^+}{B\{1 \vee \log^+ L(B)\}^{1/2}}\right)$$

$$\leq h(x_0) + h([2/\delta q]^{1/2})P(F)$$
$$+ \sup_{x \geq x_0} \frac{h(x)x^{1/2}}{\exp\{(q/2)(1-\delta)x^2\}} E\mathbb{1}\left(F^c \cap \left\{\frac{A}{B} \geq x_0\right\}\right) \frac{\exp\{(q/2)(A/B)^2\}}{(A/B)L(B)},$$

so we only require $h(x) \leq \exp(\varepsilon x^2)$ for some $\varepsilon < \frac{1}{2}$ and all large $x$ in this case. Putting $h(x) = \exp(\alpha x^2)$, with $0 < \alpha < \frac{1}{2}$, in the preceding argument then yields an absolute constant $C(\alpha)$ such that

$$(4.10) \qquad E\left[\sup_{t \geq 0} \exp\left(\frac{\alpha A_t^2}{\langle A \rangle_t \log \log(\langle A \rangle_t \vee e^2)}\right)\right] \leq C(\alpha),$$

which can be regarded as an extension to $p = 0$ of the following result of Kikuchi (1991): For every $p > 0$ and $0 < \alpha < \frac{1}{2}$, there exists an absolute constant $C_{\alpha,p}$ such that

$$E[A_\infty^{*p} \exp(\alpha A_\infty^{*2}/\langle A \rangle_\infty)] \leq C_{\alpha,p} E(A_\infty^{*p}),$$

where $A_\infty^* = \sup_{t \geq 0} |A_t|$.

By Lemma 1.5 or 1.6, (4.9) (with $t \geq 0$ replaced by $n \geq 1$) also holds for discrete-time supermartingales or martingales $A_n$ whose difference sequences satisfy the assumptions in these lemmas. Similarly, for conditionally symmetric random variables $d_i$, it follows from Lemma 1.4 and Theorem 4.1 that for every $p > 0$, there exists an absolute constant $C_p$ such that

$$(4.11) \qquad E\left(\sup_{n \geq 1} \frac{(\sum_{i=1}^n d_i)^+}{\{(\sum_{i=1}^n d_i^2) \log \log(\sum_{i=1}^n d_i^2 \vee e^2)\}^{1/2}}\right)^p \leq C_p.$$

In view of Lemma 3.9(iii), Theorem 4.1 can be applied also when $\{d_n, \mathcal{F}_n, n \geq 1\}$ is a supermartingale difference sequence such that $d_n \geq -M$ a.s. for all $n$ and some nonrandom $M > 0$. In this case, we have more generally that for $p > 0$ and $1 < r \leq 2$, there exists $C_{p,r}$ such that

$$(4.12) \quad E\left(\sup_{n \geq 1} \frac{(\sum_{i=1}^n d_i)^+}{\{(\sum_{i=1}^n |d_i|^r \vee 1)[\log \log(\sum_{i=1}^n |d_i|^r \vee e^2)]^{r-1}\}^{1/r}}\right)^p \leq C_{p,r}.$$

The remainder of this section considers maximal inequalities for self-normalized processes under condition (4.1) by using an extension of the method of mixtures introduced by Robbins and Siegmund (1970) for Brownian motion. Let $F$ be any finite measure on $(0, \lambda_0)$ with $F(0, \lambda_0) > 0$ and define the function

$$(4.13) \qquad \psi(u,v) = \int_0^{\lambda_0} \exp\{\lambda u - \lambda^r v/r\} \, dF(\lambda).$$

Given any $c > 0$ and $v > 0$, the equation $\psi(u,v) = c$ has a unique solution $u = \beta_F(v,c)$. For the case $r = 2$, the function $v \to \beta_F(v,c)$ is called a *Robbins–Siegmund boundary* in Lai (1976), in which such boundaries are shown to have the following properties:

(a) $\beta_F(v, c)$ is a concave function of $v$.
(b) $\lim_{v \to \infty} \beta_F(v, c)/v = b_F/2$, where $b_F = \sup\{y > 0 : F(0, b) = 0\}(\sup \varnothing = 0)$.
(c) If $dF(\lambda) = f(\lambda) d\lambda$ for $0 < \lambda < \lambda_0$ and $\inf_{0 < \lambda < \lambda_0} f(\lambda) > 0$ while $\sup_{0 < \lambda < \lambda_0} f(\lambda) < \infty$, then $\beta_F(v, c) \sim (v \log v)^{1/2}$ as $v \to \infty$.
(d) If $dF(\lambda) = f(\lambda) d\lambda$ for $0 < \lambda < e^{-2}$, and $= 0$ elsewhere, where

$$(4.14) \qquad f(\lambda) = 1/\{\lambda(\log \lambda^{-1})(\log\log \lambda^{-1})^{1+\delta}\},$$

for some $\delta > 0$, then as $v \to \infty$,

$$(4.15) \quad \beta_F(v, c) = \left\{2v\left[\log_2 v + \left(\frac{3}{2} + \delta\right)\log_3 v + \log\left(\frac{c}{2\sqrt{\pi}}\right) + o(1)\right]\right\}^{1/2}.$$

As in Robbins and Siegmund (1970), we write $\log_k v = \log(\log_{k-1} v)$ for $k \geq 2, \log_1 v = \log v$. For general $1 < r \leq 2$, (a) still holds, (b) holds with $b_F/2$ replaced by $b_F^{r-1}/r$ and (c) can be generalized to $\beta_F(v, c) \sim v^{1/r}\{(\log v)/(r-1)\}^{(r-1)/r}$ as $v \to \infty$. Moreover, if $f$ is given by (4.14) as in (d), then

$$(4.15) \quad \beta_F(v, c) \sim v^{1/r}\{r(\log\log v)/(r-1)\}^{(r-1)/r} \qquad \text{as } v \to \infty,$$

as can be shown by a modification of the arguments in Section 5 of Robbins and Siegmund (1970) for the case $r = 2$.

It follows from (4.1) that $\{\psi(A_t, B_t^r), t \geq 0\}$ is a nonnegative supermartingale with mean $\leq F(0, \lambda_0)$ and, therefore,

$$(4.16) \quad \begin{aligned} &P\{A_t \geq \beta_F(B_t^r, c) \text{ for some } t \geq 0\} \\ &= P\{\psi(A_t, B_t^r) \geq c \text{ for some } t \geq 0\} \leq F(0, \lambda_0)/c, \end{aligned}$$

for every $c > 0$. In particular, by choosing $c$ in (4.16) arbitrarily large, we obtain from (4.15) and (4.16) the following:

COROLLARY 4.2. *Let $1 < r \leq 2$, $\Phi_r(x) = x^r/r$ for $x \geq 0$ and suppose that (4.1) holds for the process $(A_t, B_t)$, $t \in T$. Then*

$$(4.17) \quad \limsup_{t \to \infty} \frac{A_t}{B_t(\log\log B_t)^{(r-1)/r}} \leq \left\{\frac{r}{r-1}\right\}^{(r-1)/r}$$

$$\text{a.s. on } \left\{\lim_{t \to \infty} B_t = \infty\right\}.$$

Note that Theorem 4.1 already implies the a.s. finiteness of the above lim sup on $\{\lim B_t = \infty\}$, but (4.17) gives a sharp nonrandom upper bound that reduces to the familiar $\sqrt{2}$ when $r = 2$. In view of Lemma 1.4, Corollary 4.2 with $r = 2$ is applicable to conditionally symmetric random variables $d_i$, yielding (4.17) with $A_t = \sum_{i=1}^t d_i$ and $B_t = (\sum_{i=1}^t d_i^2)^{1/2}$. The special case of this result for independent symmetric $d_i$ has been derived via an independent Rademacher sequence $\{\varepsilon_i\}$ by Griffin and Kuelbs (1991) and also by Giné and Mason (1998), who show that $\log \log B_t$ in (4.17) (with $r = 2$) can be replaced by $\log \log t$ when the $d_i$ are i.i.d. symmetric.

We next extend the preceding method of mixtures to derive maximal inequalities for conditionally symmetric $m \times 1$ vectors. An adapted sequence of random vectors $\{d_i\}$ is called *conditionally symmetric* if $\{\lambda' d_i\}$ is an adapted sequence of conditionally symmetric random variables for every $\lambda \in \mathbf{R}^m$. By Lemma 1.4, if $\{d_i\}$ is a sequence of conditionally symmetric random vectors, then for any probability distribution $F$ on $\mathbf{R}^m$, the sequence

$$(4.18) \qquad \int_{\mathbf{R}^m} \exp\left\{\lambda' \sum_{i=1}^n d_i - \tfrac{1}{2}\lambda' \sum_{i=1}^n d_i d_i' \lambda\right\} dF(\lambda), \qquad n \geq 1,$$

forms a nonnegative supermartingale with mean $\leq 1$, noting that $(\lambda' d_i)^2 = \lambda' d_i d_i' \lambda$. In particular, if we choose $F$ to be the multivariate normal distribution with mean 0 and covariance matrix $V^{-1}$, then (4.18) reduces to

$$|V|^{1/2} \left| V + \sum_{i=1}^n d_i d_i' \right|^{-1/2} \exp\left\{ \left(\sum_{i=1}^n d_i\right)' \left(V + \sum_{i=1}^n d_i d_i'\right)^{-1} \left(\sum_{i=1}^n d_i\right) \Big/ 2 \right\},$$
(4.19)

where $|\cdot|$ denotes the determinant of a square matrix. Hence, for any $c > 0$ and any positive definite $m \times m$ matrix $V$,

$$P\left\{ \frac{(\sum_{i=1}^n d_i')(V + \sum_{i=1}^n d_i d_i')^{-1}(\sum_{i=1}^n d_i)}{\log |V + \sum_{i=1}^n d_i d_i'| + 2 \log(c/\sqrt{|V|})} \geq 1 \text{ for some } n \geq 1 \right\} \leq c^{-1}.$$
(4.20)

As another application of the method of mixtures, we derive a simple formula for certain boundary crossing probabilities of multivariate continuous local martingales. Let $\lambda_{\min}(\cdot)$ denote the minimum eigenvalue of a nonnegative definite matrix.

COROLLARY 4.3. *Let $M_t$ be a continuous local martingale taking values in $\mathbf{R}^m$ such that $M_0 = 0$, $\lim_{t \to \infty} \lambda_{\min}(\langle M \rangle_t) = \infty$ a.s., and such that $E \exp(\lambda' \langle M \rangle_t \lambda) < \infty$ for all $\lambda \in \mathbf{R}^m$ and $t > 0$. Then for any $c > 1$ and any positive definite $m \times m$ matrix $V$,*

$$(4.21) \quad P\left\{ \frac{M_t'(V + \langle M \rangle_t)^{-1} M_t}{\log |V + \langle M \rangle_t| + 2 \log(c/\sqrt{|V|})} \geq 1 \text{ for some } t \geq 0 \right\} = c^{-1},$$

PROOF. First note that an expression similar to (4.19) is equal to the integral

$$(4.22) \qquad \int_{\mathbf{R}^m} \exp\left\{ \lambda' M_t - \frac{1}{2} \lambda' \langle M \rangle_t \lambda \right\} dF(\lambda),$$

where $F$ is the $m$-variate normal distribution with mean 0 and covariance matrix $V^{-1}$. Given any $\lambda \in \mathbf{R}^m$ with $\lambda \neq 0$, $\lambda' M_t$ is a univariate local martingale and $\langle \lambda' M \rangle_t = \lambda' \langle M \rangle_t \lambda \to \infty$ a.s. since $\lambda_{\min}(\langle M \rangle_t) \to \infty$ a.s. Hence, by the martingale strong law, $\lambda' M_t / \lambda' \langle M \rangle_t \lambda \to 0$ a.s. and, therefore, $\exp\{\lambda' M_t - \langle \lambda' M \rangle_t / 2\} \to 0$ a.s. as $t \to \infty$, for every $\lambda \neq 0$.

Since $E \exp(\langle \lambda' M \rangle_t / 2) < \infty$, it follows from Novikov's criterion [cf. Revuz and Yor (1999), page 332] that $\{\exp(\lambda' M_t - \langle \lambda' M \rangle_t / 2, t \geq 0\}$ is a martingale. Therefore, $\int \exp\{\lambda' M_t - \langle \lambda' M \rangle_t / 2\} dF(\lambda)$ is a nonnegative continuous martingale, and by Doob's inequality, the probability in (4.21) is $\leq c^{-1}$, similar to (4.20). Equality actually holds in (4.21), by Lemma 1 of Robbins and Siegmund (1970), if it can be shown that (4.22) converges to 0 a.s. as $t \to \infty$. Since $\exp\{\lambda' M_t - \langle \lambda' M \rangle_t\} \to 0$ a.s. for every $\lambda \neq 0$, we need only apply the dominated convergence theorem and note that by Doob's inequality,

$$P\left\{\int_{\|\lambda\| \geq a} \exp(\lambda' M_t - \langle \lambda' M \rangle_t / 2) \, dF(\lambda) \geq c \text{ for some } t \geq 0\right\}$$
$$\leq c^{-1} \int_{\|\lambda\| \geq a} dF(\lambda). \qquad \square$$

**5. A universal upper LIL.** To derive (1.7) for any adapted sequence $\{X_i\}$, one basic technique pertains to upper-bounding the probability of an event of the form $E_k = \{t_{k-1} \leq \tau_k < t_k\}$ in which $t_j$ and $\tau_j$ are stopping times defined in (5.3). Sandwiching $\tau_k$ between $t_{k-1}$ and $t_k$ enables us to replace both the random exceedance and truncation levels in (5.3) by constants. Then the event $E_k$ can be re-expressed in terms of two simultaneous inequalities, one involving centered sums and the other involving a sum of squares. Using these inequalities, we derive a supermartingale that is then used to bound $P(E_k)$. Apart from finite mean constraints, Lemma 5.1 gives the basic idea underlying the construction of this supermartingale. It will be refined in Corollary 5.3 to enable us to remove the assumptions in Lemma 5.1 concerning both the integrability of the $Y_n$'s and the restrictions on the negative part of their support.

LEMMA 5.1. *Let $0 \leq \gamma < 1$ and define*

$$(5.1) \qquad C_\gamma = -\{\gamma + \log(1-\gamma)\}/\gamma^2 = \sum_{j=2}^\infty \gamma^{j-2}\big/ j.$$

*Then $C_\gamma = c_2^{(\gamma)} = c_{\gamma,2}$, where $c_{\gamma,r}$ and $c_r^{(\gamma)}$ are the same as in Lemma 3.9. Moreover, if $Y$ is a random variable such that $Y \geq -\gamma$ and $E|Y| < \infty$, then $E \exp\{Y - EY - C_\gamma Y^2\} \leq 1$.*

PROOF. As shown in Lemma 3.9(i), $\exp(y - C_\gamma y^2) \leq 1 + y$ for all $y \geq -\gamma$. Hence, $E \exp\{Y - C_\gamma Y^2\} \leq 1 + EY \leq \exp(EY)$. $\square$

COROLLARY 5.2. *Fix any $0 \leq \gamma < 1$. Let $\{\mathcal{F}_n\}$ be an increasing sequence of $\sigma$-fields. Suppose $Y_n$ is $\mathcal{F}_n$-measurable, $E|Y_n| < \infty$ and $Y_n \geq -\gamma$ a.s. Let $\mu_n = E(Y_n | \mathcal{F}_{n-1})$. Then $\exp\{\sum_{i=1}^n (Y_i - \mu_i - C_\gamma Y_i^2)\}$ is a supermartingale whose expectation is $\leq 1$.*

COROLLARY 5.3. *Let $\{\mathcal{F}_n\}$ be an increasing sequence of $\sigma$-fields and $Y_n$ be $\mathcal{F}_n$-measurable random variables. Let $0 \leq \gamma_n < 1$ and $0 < \lambda_n \leq 1/C_{\gamma_n}$ be $\mathcal{F}_{n-1}$-measurable random variables, where $C_\gamma$ is defined in (5.1). Let $\mu_n = E\{Y_n \mathbb{1}(-\gamma_n \leq Y_n < \lambda_n) | \mathcal{F}_{n-1}\}$. Then $\exp\{\sum_{i=1}^n (Y_i - \mu_i - \lambda_i^{-1} Y_i^2)\}$ is a supermartingale whose expectation is $\leq 1$.*

PROOF. Observe that $\exp\{y - y^2/\lambda_i\} \leq 1$ if $y \geq \lambda_i$ or if $y < -\gamma_i$. Let $X_i = Y_i \mathbb{1}(-\gamma_i \leq Y_i < \lambda_i)$. Then
$$E\{\exp(Y_i - \mu_i - \lambda_i^{-1} Y_i^2) | \mathcal{F}_{i-1}\}$$
$$\leq E\{\exp(X_i - \mu_i - \lambda_i^{-1} X_i^2) | \mathcal{F}_{i-1}\}$$
$$\leq E\{(1 + X_i) e^{-\mu_i} | \mathcal{F}_{i-1}\} = (1 + \mu_i) e^{-\mu_i};$$
see the proof of Lemma 5.1 for the last inequality, recalling that $\mu_i = E(X_i | \mathcal{F}_{i-1})$. Since $(1 + x) e^{-x} \leq 1$ for all $x$, the desired conclusion follows. □

The centering constants in (1.7) involve sums of expectations conditioned on the past which are computed as functions of the endpoints of the interval on which the associated random variable is truncated. The actual endpoints used, however, are neither knowable nor determined until the future. Thus the centered sums that result are not a martingale. Nevertheless, by using certain stopping times, the random truncation levels can be replaced by non-random ones, thereby yielding a supermartingale structure for which Corollary 5.5 applies, enabling us to establish the following result.

THEOREM 5.4. *Let $X_n$ be measurable with respect to $\mathcal{F}_n$, an increasing sequence of $\sigma$-fields. Let $\lambda > 0$ and $h(\lambda)$ be the positive solution of*

(5.2) $$h - \log(1 + h) = \lambda^2.$$

*Let $b_\lambda = h(\lambda)/\lambda$, $\gamma = h(\lambda)/\{1 + h(\lambda)\}$ and $a_\lambda = \lambda/(\gamma C_\gamma)$, where $C_\gamma$ is defined by (5.1). Then (1.7) holds on $\{\lim_{n \to \infty} V_n = \infty\}$ and $\lim_{\lambda \to 0} b_\lambda = \sqrt{2}$.*

PROOF. Recall that $V_n^2 = X_1^2 + \cdots + X_n^2$ and $v_n = V_n (\log \log V_n)^{-1/2}$. Let $e_k = \exp(k/\log k)$. Define

$$t_j = \inf\{n : V_n \geq e_j\},$$

(5.3) $$\tau_j = \inf\left\{n \geq t_j : S_n - \sum_{i=1}^n \mu_i(-\lambda v_n, a_\lambda v_n) \geq (1 + 3\varepsilon) b_\lambda V_n (\log \log V_n)^{1/2}\right\},$$

letting $\inf \varnothing = \infty$. To prove (1.7), it suffices to show that for all sufficiently small $\varepsilon > 0$,

$$\lim_{K \to \infty} \sum_{k=K}^{\infty} P\{\tau_k < t_{k+1}\} = 0. \tag{5.4}$$

Note that $\tau_k \geq t_k$ and that $t_k$ may equal $t_{k+1}$, in which case $\{\tau_k < t_{k+1}\}$ becomes the empty set. Moreover, on $\{\lim_{n\to\infty} V_n = \infty\}$, $t_j < \infty$ for every $j$ and $\lim_{j\to\infty} t_j = \infty$. Since $y(\log \log y)^{-1/2}$ is increasing in $y \geq e_3$, we have the following inequalities on $\{t_k \leq \tau_k < t_{k+1}\}$ with $k \geq 3$:

$$e_k \leq \left( \sum_{i=1}^{\tau_k} X_i^2 \right)^{1/2} < e_{k+1}, \tag{5.5}$$

$$d_k := e_k (\log \log e_k)^{-1/2} \leq v_{t_k} \leq v_{\tau_k} < d_{k+1}, \tag{5.6}$$

$$\mu_i(-\lambda v_{\tau_k}, a_\lambda v_{\tau_k}) \geq \mu_i(-\lambda d_{k+1}, a_\lambda d_k) \quad \text{for } 1 \leq i \leq \tau_k. \tag{5.7}$$

Let $\mu_{i,k} = \mu_i(-\lambda d_{k+1}, a_\lambda d_k)$. We shall replace $X_i$ (for $1 \leq i \leq \tau_k$) by $Y_{i,k} := (\lambda d_{k+1})^{-1} \gamma X_i$ and $\mu_{i,k}$ by $\tilde{\mu}_{i,k} := (\lambda d_{k+1})^{-1} \gamma \mu_i(-\lambda d_{k+1}, a_\lambda d_k)$. Since $\lambda^{-1} \gamma a_\lambda = C_\gamma^{-1}$,

$$\tilde{\mu}_{i,k} = E\{Y_{i,k} \mathbb{1}(-\gamma \leq Y_{i,k} < C_\gamma^{-1} d_k / d_{k+1}) | \mathcal{F}_{i-1}\}. \tag{5.8}$$

Since $e_k / d_k = (\log \log e_k)^{1/2}$ and $d_k / d_{k+1} \to 1$ as $k \to \infty$, it follows from (5.5)–(5.7) that for all sufficiently large $k$, the event $\{t_k \leq \tau_k < t_{k+1}\}$ is a subset of

$$\left\{ \sum_{i=1}^{\tau_k} (\lambda d_{k+1})^{-1}(X_i - \mu_{i,k}) \geq (1 + 2\varepsilon)\lambda^{-1} b_\lambda \log \log e_k, \ \tau_k < \infty \right\}$$

$$\subset \left\{ \sum_{i=1}^{\tau_k} [(\lambda d_{k+1})^{-1} \gamma (X_i - \mu_{i,k}) - C_\gamma (d_{k+1}/d_k)(\lambda d_{k+1})^{-2} \gamma^2 X_i^2] \right.$$

$$\left. \geq (1 + 2\varepsilon)\gamma \lambda^{-1} b_\lambda \log \log e_k - C_\gamma (d_{k+1}/d_k)(\gamma/\lambda)^2 \log \log e_{k+1}, \tau_k < \infty \right\}$$

$$\subset \left\{ \sup_{n \geq 1} \exp\left[ \sum_{i=1}^{n} (Y_{i,k} - \tilde{\mu}_{i,k} - C_\gamma d_k^{-1} d_{k+1} Y_{i,k}^2) \right] \right.$$

$$\left. \geq \exp[(1 + \varepsilon)(\gamma \lambda^{-1} b_\lambda - C_\gamma \gamma^2 \lambda^{-2})(\log k)] \right\}.$$

In view of (5.8), we can apply Corollary 5.3 to conclude that the last event above involves the supremum of a nonnegative supermartingale with mean $\leq 1$. Therefore, application of Doob's inequality to this event yields

$$P\{\tau_k < t_{k+1}\} \leq \exp\{-(1 + \varepsilon)(\gamma \lambda^{-1} b_\lambda - C_\gamma \gamma^2 \lambda^{-2})(\log k)\},$$

which implies (5.4) since

(5.9) $\quad \gamma\lambda^{-1}b_\lambda - \lambda^{-2}\gamma^2 C_\gamma = \lambda^{-2}\{\gamma h(\lambda) + \gamma + \log(1-\gamma)\} = 1.$

The first equality in (5.9) follows from (5.1) and $b_\lambda = h(\lambda)/\lambda$, and the second equality from $\gamma = h(\lambda)/(1 + h(\lambda))$ and (5.2). Moreover, (5.2) implies that $h^2(\lambda) \sim 2\lambda^2$ and, therefore, $b_\lambda \to \sqrt{2}$ as $\lambda \to 0$. □

REMARK 5.5. The choice of $\gamma$ in Theorem 5.4 actually comes from minimizing $\gamma\lambda^{-1}b_\lambda - \lambda^{-2}\gamma^2 C_\gamma$ over $0 < \gamma < 1$, whereas $b_\lambda$ is employed to make this minimizing value equal to 1, leading to the equation (5.2) defining $h(\lambda)$.

As pointed out in Section 1, an immediate consequence of Theorem 5.4 is the upper half (1.7) of the LIL for any supermartingale whose difference sequence $X_n$ is bounded below by $-\lambda v_n$. The following example shows that we cannot dispense with this boundedness assumption.

EXAMPLE 5.6. Let $X_1 = X_2 = 0, X_3, X_4, \ldots$ be independent random variables such that

$P\{X_n = -n^{-1/2}\} = 1/2 - n^{-1/2}(\log n)^{1/2} - n^{-1}(\log n)^{-2},$

$P\{X_n = -m_n\} = n^{-1}(\log n)^{-2}, \qquad P\{X_n = n^{-1/2}\} = 1/2 + n^{-1/2}(\log n)^{1/2}$

for $n \geq 3$, where $m_n \sim 2(\log n)^{5/2}$ is chosen so that $EX_n = 0$. Then $P\{X_n = -m_n \text{ i.o.}\} = 0$. Hence, with probability $1, V_n^2 = \sum_{i=1}^n i^{-1} + O(1) = \log n + O(1)$. Since $\widetilde{X}_i := X_i \mathbb{1}(|X_i| \leq 1) - EX_i \mathbb{1}(|X_i| \leq 1)$ are independent bounded random variables with zero means and $\text{Var}(\widetilde{X}_i) \sim i^{-1}$, Kolmogorov's LIL yields

(5.10) $\quad \limsup_{n\to\infty}\left(\sum_{i=1}^n \widetilde{X}_i\right)\bigg/\{2(\log n)(\log\log\log n)\}^{1/2} = 1 \quad \text{a.s.}$

Since $\sum_{i=1}^n EX_i \mathbb{1}(|X_i| \leq 1) \sim 2\sum_{i=1}^n i^{-1}(\log i)^{1/2} \sim \frac{4}{3}(\log n)^{3/2}$, this implies that with probability 1,

$$\frac{\sum_{i=1}^n X_i}{V_n (\log\log V_n)^{1/2}} \sim \frac{\sum_{i=1}^n X_i \mathbb{1}(|X_i| \leq 1)}{\{(\log n)(\log\log\log n)\}^{1/2}}$$

$$\sim \frac{4(\log n)^{3/2}}{3\{(\log n)(\log\log\log n)\}^{1/2}} \to \infty.$$

Note that $m_n(\log\log V_n)^{1/2}/V_n \to \infty$. This shows that without the boundedness condition $X_n \geq -\lambda V_n(\log\log V_n)^{-1/2}$, the upper LIL need not hold for martingales self-normalized by $V_n$. It also shows the importance of the centering in Theorem 5.4 because subtracting $EX_i \mathbb{1}(|X_i| \leq 1)$ from $X_i$ gives the LIL in view of (5.10).

Note that Corollary 5.3, which leads to Theorem 5.4, only uses the special case $r = 2$ of Lemma 3.9(i). More generally, for $1 < r \leq 2$, we can use Lemma 3.9(i) and the same arguments as those in Lemma 5.1 and Corollary 5.3 to show that

$$\exp\left\{\sum_{i=1}^{n}(Y_i - E[Y_i \mathbb{1}(-\gamma_i \leq Y_i < \lambda_i^{1/(r-1)})|\mathcal{F}_{i-1}] - \lambda_i^{-1}|Y_i|^r)\right\}, \tag{5.11}$$

$n \geq 1$, is a supermartingale,

for any $\mathcal{F}_{i-1}$-measurable random variables $0 \leq \gamma_i < 1$ and $0 < \lambda_i \leq 1/c_{\gamma,r}$, where $c_{\gamma,r}$ is defined in Lemma 3.9. Therefore, Theorem 5.4 can be extended to the following:

THEOREM 5.7. *Let $X_n$ be measurable with respect to $\mathcal{F}_n$, an increasing sequence of $\sigma$-fields. For $1 < r \leq 2$, let $V_{n,r} = (\sum_{i=1}^{n}|X_i|^r)^{1/r}$, $v_{n,r} = V_{n,r}\{\log\log(V_{n,r} \vee e^2)\}^{-1/r}$. Then for any $0 < \gamma < 1$, there exists a positive constant $b_{\gamma,r}$ such that*

$$\limsup_{n \to \infty}\left\{S_n - \sum_{i=1}^{n}\mu_i(-\gamma v_{n,r}, c_{\gamma,r}^{-1/(r-1)}v_{n,r})\right\}\Big/\{V_{n,r}(\log\log V_{n,r})^{(r-1)/r}\}$$

$$\leq b_{\gamma,r} \quad a.s.$$

*on $\{\lim_{n \to \infty} V_{n,r} = \infty\}$, where $c_{\gamma,r}$ is given in Lemma 3.9.*

**6. Compact LIL for self-normalized martingales and applications to sums of independent random variables.** Although Theorem 5.4 gives an upper LIL for any adapted sequence $\{X_i\}$, the upper bound in (1.7) may not be attained. A simple example is $S_n = \sum_{i=1}^{n} w_i Y_i$, where $w_i = i!$ and $Y_1, Y_2, \ldots$ are i.i.d. with $P\{Y_i = 1\} = \frac{1}{2} = P\{Y_i = -1\}$. Here $V_n = (\sum_{i=1}^{n} w_i^2 Y_i^2)^{1/2} \sim n!$, $S_n/V_n = \mathrm{sgn}(Y_n) + o(1)$ and $\sum_{i=1}^{n}\mu_i(-\lambda v_n, a_\lambda v_n) = o(V_n)$ a.s. Thus, the norming term $V_n(\log\log V_n)^{1/2}$ is too large in this case. In this section we consider the case of martingales $\{S_n, \mathcal{F}_n, n \geq 1\}$ self-normalized by $V_n$ and prove the lower half counterpart of (1.8) when the increments of $S_n$ do not grow too fast, thereby establishing (1.9). This is the content of Theorem 6.1, which is further strengthened into a compact LIL in Corollary 6.2. We end this section with an application to weighted sums (with random weights) of i.i.d. random variables, a remark on Theorem 5.4 and an example highlighting the difference between this LIL and an analogous LIL of Stout (1970) in which $V_n^2$ is replaced by $s_n^2 = \sum_{i=1}^{n} E(X_i^2|\mathcal{F}_{i-1})$.

THEOREM 6.1. *Let $\{X_n\}$ be a martingale difference sequence with respect to an increasing sequence of $\sigma$-fields $\mathcal{F}_n$ such that $|X_n| \leq m_n$ a.s. for some $\mathcal{F}_{n-1}$-measurable random variable $m_n$, with $V_n \to \infty$ and $m_n/\{V_n(\log\log V_n)^{-1/2}\} \to 0$ a.s. Then (1.9) holds.*

PROOF. Take $0 < b < \beta < \tilde{\beta} < \sqrt{2}$. Since $1 - \Phi(x) = \exp\{-(\frac{1}{2} + o(1))x^2\}$ as $x \to \infty$, we can choose $\lambda$ sufficiently large such that

$$(6.1) \qquad \{1 - \Phi(\beta\sqrt{\lambda})\}^{1/\lambda} \geq \exp(-\tilde{\beta}^2/2),$$

where $\Phi$ is the standard normal distribution function. Take $a > 1$ and define for $j \geq 2$ and $k = 0, 1, \ldots, [\lambda^{-1} \log j]$,

$$a_{j,k} = a^j + k(a^{j+1} - a^j)/[\lambda^{-1} \log j], \qquad t_j(k) = \inf\{n : V_n^2 \geq a_{j,k}\}.$$

Let $t_j = \inf\{n : V_n^2 \geq a^j\}$, so $t_j(0) = t_j$, $t_j([\lambda^{-1} \log j]) = t_{j+1}$. Since $X_n^2 = o(V_n^2 (\log \log V_n)^{-1})$ a.s. and $a_{j,k} \leq V_{t_j(k)}^2 < a_{j,k} + X_{t_j(k)}^2$,

$$(6.2) \qquad V_{t_j(k)}^2 = a_{j,k}\{1 + o((\log j)^{-1})\} \qquad \text{a.s.}$$

It will be shown that

$$(6.3) \qquad \sum_{t_j(k) < n \leq t_j(k+1)} X_n^2 \Big/ \sum_{t_j(k) < n \leq t_j(k+1)} E(X_n^2 | \mathcal{F}_{n-1}) \to 1$$

$$\text{in probability under } P(\cdot | \mathcal{F}_{t_j(k)})$$

as $j \to \infty$, uniformly in $0 \leq k < [\lambda^{-1} \log j]$.

Let $S_{m,n} = \sum_{m < i \leq n} X_i$, $V_{m,n}^2 = \sum_{m < i \leq n} X_i^2$. In view of (6.2),

$$(6.4) \qquad V_{t_j(k), t_j(k+1)}^2 \sim a^j(a-1)/[\lambda^{-1} \log j], \qquad V_{t_j, t_{j+1}}^2 \sim a^j(a-1) \quad \text{a.s.}$$

Since $X_n^2$ is bounded by the $\mathcal{F}_{n-1}$-measurable random variable $m_n^2$, which is $o(V_n^2 (\log \log V_n)^{-1})$ a.s., the conditional Lindeberg condition holds and, in view of (6.3) and (6.4), the martingale central limit theorem [cf. Durrett (1996), page 414] can be applied to yield

$$(6.5) \quad P\{S_{t_j(k), t_j(k+1)} \geq \beta\sqrt{\lambda} V_{t_j(k), t_j(k+1)} | \mathcal{F}_{t_j(k)}\} \to 1 - \Phi(\beta\sqrt{\lambda}) \qquad \text{a.s.}$$

as $j \to \infty$, uniformly in $0 \leq k < [\lambda^{-1} \log j]$. Since

$$S_{t_j, t_{j+1}} = \sum_{0 \leq k < [\lambda^{-1} \log j]} S_{t_j(k), t_j(k+1)}$$

and

$$V_{t_j, t_{j+1}} (\log j)^{1/2} = (\sqrt{\lambda} + o(1)) \sum_{0 \leq k < [\lambda^{-1} \log j]} V_{t_j(k), t_j(k+1)} \qquad \text{a.s.}$$

by (6.4), it follows from (6.5) that as $j \to \infty$,

$$P\{S_{t_j, t_{j+1}} \geq b V_{t_j, t_{j+1}} (\log j)^{1/2} | \mathcal{F}_{t_j}\}$$
$$\geq P\{S_{t_j(k), t_j(k+1)} \geq \beta\sqrt{\lambda} V_{t_j(k), t_j(k+1)} \text{ for all } 0 \leq k < [\lambda^{-1} \log j] | \mathcal{F}_{t_j}\}$$
$$= (1 - \Phi(\beta\sqrt{\lambda}) + o(1))^{[\lambda^{-1} \log j]}$$
$$\geq \exp\{-(\tilde{\beta}^2/2 + o(1)) \log j\} \qquad \text{a.s.,}$$

in view of (6.1). Since $\tilde{\beta}^2/2 < 1$, the conditional Borel–Cantelli lemma then yields

(6.6) $$\limsup_{j \to \infty} S_{t_j, t_{j+1}} \Big/ \{V_{t_j, t_{j+1}} (\log j)^{1/2}\} \geq b \quad \text{a.s.}$$

Recalling that $V_n \to \infty$ and $m_n = o(V_n (\log \log V_n)^{-1/2})$ a.s., we obtain from (1.8) that

(6.7) $$\limsup_{n \to \infty} S_n \Big/ \{V_n (\log \log V_n)^{1/2}\} \leq \sqrt{2} \quad \text{a.s.},$$

and the same conclusion still holds with $S_n$ replaced by $-S_n$ (which is a martingale). Combining this with (6.4) and (6.6) yields

(6.8) $$\limsup_{j \to \infty} S_{t_{j+1}} \Big/ \{V_{t_{j+1}} (\log \log V_{t_{j+1}})^{1/2}\}$$
$$\geq b a^{-1/2} (a-1)^{1/2} - \sqrt{2} a^{-1/2} \quad \text{a.s.}$$

Since $a$ can be chosen arbitrarily large and $b$ arbitrarily close to $\sqrt{2}$ in (6.8),

$$\limsup_{j \to \infty} S_{t_{j+1}} \Big/ \{V_{t_{j+1}} (\log \log V_{t_{j+1}})^{1/2}\} \geq \sqrt{2} \quad \text{a.s.}$$

Combining this with the upper half result (6.7) yields (1.9).

It remains to prove (6.3). Let $\alpha_j = a^j(a-1)/[\lambda^{-1} \log j]$. In view of (6.4), we need to show that given any $0 < \rho < \frac{1}{2}$ and $\delta > 0$,

(6.9) $$\limsup \left[ P\left\{ \sum_{t_j(k) < n \leq t_j(k+1)} E(X_n^2 | \mathcal{F}_{n-1}) \geq (1+\rho)\alpha_j \Big| \mathcal{F}_{t_j(k)} \right\} \right.$$
$$\left. + P\left\{ \sum_{t_j(k) < n \leq t_j(k+1)} E(X_n^2 | \mathcal{F}_{n-1}) \leq (1-\rho)\alpha_j \Big| \mathcal{F}_{t_j(k)} \right\} \right] \leq \delta \quad \text{a.s.}$$

Choose $\varepsilon > 0$ such that $2\{\max[(1+\rho)e^{-\rho}, (1-\rho)e^\rho]\}^{1/\varepsilon} < \delta$. Let $\tilde{X}_n = X_n \mathbb{1}(m_n^2 \leq \varepsilon \alpha_j)$ and note that since $m_n$ is $\mathcal{F}_{n-1}$-measurable and $X_n^2 \leq m_n^2$,

$$0 \leq E(X_n^2 | \mathcal{F}_{n-1}) - E(\tilde{X}_n^2 | \mathcal{F}_{n-1}) \leq m_n^2 \mathbb{1}(m_n^2 > \varepsilon \alpha_j).$$

Moreover, $P\{m_n^2 \leq \varepsilon \alpha_j \text{ for all } t_j(k) < n \leq t_j(k+1) | \mathcal{F}_{t_j(k)}\} \to 1$ a.s. Hence, it suffices to consider $E(\tilde{X}_n^2 | \mathcal{F}_{n-1})$ instead of $E(X_n^2 | \mathcal{F}_{n-1})$ in (6.9). Since $\tilde{X}_n^2 \leq \varepsilon \alpha_j$, we can apply Corollary 15 of Freedman (1973) to conclude that

$$P\left\{ \sum_{t_j(k) < n \leq t_j(k+1)} E(\tilde{X}_n^2 | \mathcal{F}_{n-1}) \geq (1+\rho)\alpha_j \Big| \mathcal{F}_{t_j(k)} \right\}$$
$$+ P\left\{ \sum_{t_j(k) < n \leq t_j(k+1)} E(\tilde{X}_n^2 | \mathcal{F}_{n-1}) \leq (1-\rho)\alpha_j \Big| \mathcal{F}_{t_j(k)} \right\}$$
$$\leq (1+\rho)e^{-\rho/\varepsilon} + (1-\rho)e^{\rho/\varepsilon} + o(1) < \delta,$$

completing the proof. □

COROLLARY 6.2. *With the same notation and assumptions as in Theorem 6.1, the cluster set of the sequence $\{S_n/[V_n(\log\log(V_n \vee e^2))^{1/2}]\}$ is the interval $[-\sqrt{2}, \sqrt{2}]$.*

PROOF. Replacing $X_n$ by $-X_n$ in Theorem 6.1 yields $\liminf_{n\to\infty} S_n/\{V_n(\log\log V_n)^{1/2}\} = -\sqrt{2}$ a.s. The desired conclusion then follows from Proposition 2.1 of Griffin and Kuelbs (1989). □

EXAMPLE 6.3. Let $Y_1, Y_2, \ldots$ be i.i.d. random variables with a common distribution function $F$ having mean 0. Let $\mathcal{F}_n$ be the $\sigma$-field generated by $Y_1, \ldots, Y_n$. Let $w_n$ be $\mathcal{F}_{n-1}$-measurable and let $S_n = \sum_{i=1}^n w_i Y_i$, $V_n^2 = \sum_{i=1}^n w_i^2 Y_i^2$. Suppose $V_n \to \infty$ a.s. and there exists $\mathcal{F}_{n-1}$-measurable $m_n$ such that with probability 1,

$$(6.10) \qquad 0 < m_n = o(V_n(\log\log V_n)^{-1/2}),$$

$$(6.11) \qquad \sum_{i=1}^n w_i \int_{|w_i x| \geq m_i} x\, dF(x) = o(V_n(\log\log V_n)^{1/2}),$$

$$\sum_{i=1}^n \left\{ w_i \int_{|w_i x| \geq m_i} x\, dF(x) \right\}^2 = o(V_n^2),$$

$$(6.12) \qquad \sum_{n=1}^\infty \left\{ \bar{F}(m_n/|w_n|) + F(-m_n/|w_n|) \right\} < \infty,$$

where $\bar{F}(x) = P(Y_i \geq x) = 1 - F(x-)$. Let $X_n = w_n Y_n \mathbb{1}(|w_n Y_n| < m_n)$. Then $E(X_n|\mathcal{F}_{n-1}) = -w_n \int_{|w_n x| \geq m_n} x\, dF(x)$. Moreover, by (6.12) and the conditional Borel–Cantelli lemma, with probability 1,

$$(6.13) \qquad w_n Y_n = X_n \text{ for all large } n \text{ and therefore } V_n^2 = \sum_{i=1}^n X_i^2 + O(1).$$

Applying Corollary 6.2 to $\sum_{i=1}^n \{X_i - E(X_i|\mathcal{F}_{i-1})\}$ (with $|X_i| < m_i$) and combining the result with (6.11) and (6.13), we obtain $[-\sqrt{2}, \sqrt{2}]$ as the a.s. cluster set of the sequence $\{S_n/[V_n(\log\log(V_n \vee e^2))^{1/2}]\}$. Note in this connection that

$$\sum_{i=1}^n \{X_i - E(X_i|\mathcal{F}_{i-1})\}^2$$

$$= \sum_{i=1}^n X_i^2 - 2\sum_{i=1}^n \{X_i - E(X_i|\mathcal{F}_{i-1})\}E(X_i|\mathcal{F}_{i-1}) - \sum_{i=1}^n E^2(X_i|\mathcal{F}_{i-1})$$

$$= \sum_{i=1}^{n} X_i^2 - \sum_{i=1}^{n} E^2(X_i|\mathcal{F}_{i-1})$$

$$+ O\left(\left(\sum_{i=1}^{n}\{X_i - E(X_i|\mathcal{F}_{i-1})\}^2\right)^{1/2}\left(\sum_{i=1}^{n} E^2(X_i|\mathcal{F}_{i-1})\right)^{1/2}\right).$$

Note that Theorems 6.1 and 6.2 pertain to martingale difference sequences $X_n$. This means that given an integrable sequence $\{X_n\}$, one should first consider centering $X_n$ at its conditional expectation given $\mathcal{F}_{n-1}$ before applying the theorems to $\widetilde{X}_n = X_n - E(X_n|\mathcal{F}_{n-1})$ and $V_n = (\sum_{i=1}^{n} \widetilde{X}_i^2)^{1/2}$. Although Theorem 6.1 requires $\widetilde{X}_n$ to be bounded by $\mathcal{F}_{n-1}$-measurable $m_n = o(V_n(\log \log V_n)^{-1/2})$, we can often dispense with such boundedness assumption via a truncation argument, as shown in Example 6.3. In the more general context of Theorem 5.4, the $X_n$ may not be even integrable, so Theorem 5.4 centers the $X_n$ at certain truncated conditional expectations. Using $(\sum_{i=1}^{n} X_i^2)^{1/2}$ for the norming factor, however, may be too large since it involves uncentered $X_i$'s. To alleviate this problem, we can first center $X_n$ at its conditional median before applying Theorem 5.4 to $\widetilde{X}_n = X_n - \text{med}(X_n|\mathcal{F}_{n-1})$, as illustrated in the following:

EXAMPLE 6.4. Let $0 < \alpha < 1$, $d_1 \geq 0$, $d_2 \geq 0$ with $d_1 + d_2 > 0$. Let $Y, Y_1, Y_2, \ldots$ be i.i.d. random variables such that

(6.14)
$$P\{Y \geq y\} = (d_1 + o(1))y^{-\alpha},$$
$$P\{Y \leq -y\} = (d_2 + o(1))y^{-\alpha} \qquad \text{as } y \to \infty.$$

Let $\widehat{S}_n = \sum_{i=1}^{n} Y_i$, $\widehat{V}_n^2 = \sum_{i=1}^{n} Y_i^2$, $\widehat{v}_n = \widehat{V}_n(\log\log \widehat{V}_n)^{-1/2}$. Then by Theorem 5.1 of Shao (1997),

(6.15) $\quad \limsup_{n\to\infty} \widehat{S}_n/\{\widehat{V}_n(\log\log n)^{1/2}\} = \{\beta(\alpha, d_1, d_2)\}^{-1/2}$ a.s.

for some positive constant $\beta(\alpha, d_1, d_2)$ which is given explicitly in his Theorem 3.2. Moreover, $E\{Y\mathbb{1}(-\lambda y \leq Y < a_\lambda y)\} = (d_1 a_\lambda - d_2\lambda + o(1))\alpha y^{1-\alpha}/(1-\alpha)$ as $y \to \infty$ and

(6.16)
$$n\widehat{v}_n^{1-\alpha}/\{\widehat{V}_n(\log\log \widehat{V}_n)^{1/2}\} = n/\{\widehat{V}_n^\alpha(\log\log \widehat{V}_n)^{(2-\alpha)/2}\} = O(1) \qquad \text{a.s.}$$

since $\log\log \widehat{V}_n \sim \log\log n$ and

$$\liminf_{n\to\infty}\left(\sum_{i=1}^{n} Y_i^2\right)\bigg/\{n^{1/\widetilde{\alpha}}(\log\log n)^{-(1-\widetilde{\alpha})/\widetilde{\alpha}}\} > 0 \qquad \text{a.s. with } \widetilde{\alpha} = \alpha/2,$$

by the so-called delicate LIL [cf. Breiman (1968)].

Now let $X_n = n^r + Y_n$ with $r > 1/\alpha$ and let $S_n = \sum_{i=1}^{n} X_i$, $V_n^2 = \sum_{i=1}^{n} X_i^2$. Since $Y_n = o(n^s)$ a.s. for any $s > 1/\alpha$, it follows that $S_n \sim V_n \sim n^{r+1}/(r+$

1) and $\mu_i(-\lambda v_n, a_\lambda v_n) = i^r + o(n^{(r+1)(1-\alpha)}) = i^r + o(n^r)$ a.s., recalling that $r\alpha > 1$. Therefore, although (1.7) still holds in this case, it is too crude as the nonrandom location shift $n^r$ is the dominant term in $X_n$ causing $V_n$ to swamp the centered $S_n$. Centering the $X_n$ first at its median will remove this problem. Specifically, if we apply (1.7) to $\widetilde{X}_n = X_n - \mathrm{med}(X_n)$ and $\widetilde{V}_n^2 = \sum_{i=1}^n \widetilde{X}_i^2$, then $\widetilde{X}_n = Y_n - \mathrm{med}(Y)$ and (6.15) still holds with $\widehat{S}_n$ replaced by $\widetilde{S}_n$.

The following example shows that one cannot dispense with the assumptions of Theorem 6.1 and highlights the difference between our result and the LIL of Stout (1970), where the martingale $S_n$ is normalized by the square root of the conditional variance $\sum_{i=1}^n E(X_i^2 | \mathcal{F}_{i-1})$.

EXAMPLE 6.5. Taking $X_1 = 0, X_2, X_3, \ldots$ and $m_n$ as in Example 5.5, let $Y_n = X_n \mathbb{1}(|X_n| \leq 1)$. Then $P\{Y_n \neq X_n \text{ i.o.}\} = P\{Y_n = -m_n \text{ i.o.}\} = 0$. As shown in Example 5.5, with probability 1, $V_n^2 = \sum_{i=1}^n Y_i^2 + O(1) = \log n + O(1)$ and

$$(6.17) \qquad \frac{\sum_{i=1}^n X_i}{V_n (\log \log V_n)^{1/2}} \sim \frac{4(\log n)^{3/2}}{3\{(\log n)(\log \log \log n)\}^{1/2}} \to \infty.$$

Note that $m_n (\log \log V_n)^{1/2}/V_n \to \infty$. This shows that without the condition $m_n/\{V_n(\log \log V_n)^{-1/2}\} \to 0$, the LIL need not hold for martingales self-normalized by $V_n$. On the other hand, $X_n$ is clearly bounded above and, therefore, satisfies the boundedness condition of Stout (1970). Note that $\mathrm{Var}(X_i) \sim 4(\log i)^3/i$ and, therefore, $s_n^2 := \sum_{i=1}^n E(X_i^2|\mathcal{F}_{i-1}) \sim (\log n)^4$, yielding

$$(6.18) \qquad \frac{\sum_{i=1}^n X_i}{s_n(\log \log s_n)^{1/2}} \sim \frac{4(\log n)^{3/2}}{3(\log n)^2 (\log \log \log n)^{1/2}} \to 0 \qquad \text{a.s.,}$$

which is consistent with Stout's (1970) upper LIL. Contrasting (6.18) with (6.17) shows the difference between Stout's result and ours. Notice that what is being investigated in (6.17) is the maximal a.s. growth rate of $S_n$. To assess it we employed a norming sequence based on the square root of its sum of squares. This technique works properly only when $S_n$ is adequately centered, as in (1.7). By contrast, in the approach of Stout, a norming sequence is generated from the square root of the sum of conditional expectations of these squares. However, in the absence of a suitable truncation of the random variables this quantity is also inappropriate for investigating almost sure behavior whenever expectations overinflate the impact of large values of the squares which occur too infrequently to be relevant with respect to almost sure behavior.

V. H. de la Peña
Department of Statistics
Columbia University
New York, New York 10027
USA
e-mail: vhdl@columbia.edu

M. J. Klass
Departments of Statistics
and Mathematics
University of California at Berkeley
Berkeley, California 94720
USA

T. L. Lai
Department of Statistics
Stanford University
Stanford, California 94305
USA
e-mail: lait@stat.stanford.edu